\documentclass[journal]{IEEEtran}
\usepackage{multirow}
\usepackage{multicol}
\usepackage{lipsum}
\usepackage{graphics}
\usepackage{amsmath}
\usepackage{amsbsy}
\usepackage{blindtext}
\usepackage{tcolorbox}
\usepackage{amssymb}
\usepackage{scalerel}
\usepackage{mathabx}
\usepackage{verbatim}
\usepackage{booktabs}
\usepackage{rotating,tabularx}
\usepackage{mathrsfs} 
\usepackage{hyperref}

\usepackage{soul}
\usepackage{xcolor}
\usepackage{cite}
\usepackage{tikz}
\usepackage{color, colortbl}
\usepackage[first=0,last=9]{lcg}
\definecolor{LightGray}{rgb}{0.7,0.7,0.7}

\usepackage[wby]{callouts}
\usetikzlibrary{shapes.multipart}

\usepackage{array}
\usepackage{makecell}

\allowdisplaybreaks

\usepackage{amsthm}
\theoremstyle{definition}

\theoremstyle{remark}

\usepackage[utf8]{inputenc}
\usepackage[english]{babel}
 
\newtheorem{proposition}{Proposition}

\usepackage[caption=false,font=footnotesize]{subfig}
\usepackage[T1]{fontenc}
\usepackage{scalerel,stackengine}
\newcommand\reallywidecheck[1]{%
\savestack{\tmpbox}{\stretchto{%
  \scaleto{%
    \scalerel*[\widthof{\ensuremath{#1}}]{\kern-.6pt\bigwedge\kern-.6pt}%
    {\rule[-\textheight/2]{1ex}{\textheight}}
  }{\textheight}%
}{0.5ex}}%
\stackon[1pt]{#1}{\scalebox{-1}{\tmpbox}}%
}
\stackMath
\IEEEoverridecommandlockouts
\IEEEoverridecommandlockouts
\renewcommand\Re{\operatorname{Re}}
\renewcommand\Im{\operatorname{Im}}

\newif\ifarxiv
 \arxivtrue

\ifarxiv
\newcommand*{\mrn}{\textcolor{black}}

\else
\newcommand*{\mrn}{\textcolor{black}}

\fi

\title{Tightening QC Relaxations of AC Optimal Power Flow Problems via Complex Per Unit Normalization}

\ifarxiv
\author{Mohammad Rasoul Narimani,$^{\ast}$ Daniel K. Molzahn,$^{\dagger}$ and Mariesa L. Crow$^{\ast}$%
\thanks{${\ast}$: Electrical and Computer Engineering Department, Missouri University of Science and Technology. \{mn9t5, crow\}@mst.edu.}%
\thanks{$^{\ddagger}$: School of Electrical and Computer Engineering, Georgia Institute of Technology. molzahn@gatech.edu.}}
\else
\author{Mohammad Rasoul Narimani, Daniel K. Molzahn, and Mariesa L. Crow}%
\fi

\begin{document}
\maketitle

\begin{abstract}
Optimal power flow (OPF) is a key problem in power system operations. OPF problems that use the nonlinear AC power flow equations to accurately model the network physics have inherent challenges associated with non-convexity. To address these challenges, recent research has applied various convex relaxation approaches to OPF problems. The QC relaxation is a promising approach that convexifies the trigonometric and product terms in the OPF problem by enclosing these terms in convex envelopes. The accuracy of the QC relaxation strongly depends on the tightness of these envelopes. This paper presents two improvements to these envelopes. The first improvement leverages a polar representation of the branch admittances in addition to the rectangular representation used previously. The second improvement is based on a coordinate transformation via a complex per unit base power normalization that rotates the power flow equations. The trigonometric envelopes resulting from this rotation can be tighter than the corresponding envelopes in previous QC relaxation formulations. Using an empirical analysis with a variety of test cases, this paper suggests an appropriate value for the angle of the complex base power. Comparing the results with a state-of-the-art QC formulation reveals the advantages of the proposed improvements.
\end{abstract}

\begin{IEEEkeywords}
Optimal power flow, Convex relaxation
\end{IEEEkeywords}

\IEEEpeerreviewmaketitle

\section{Introduction}

\IEEEPARstart{O}{ptimal} power flow (OPF) problems are central to many tasks in power system operations. OPF problems optimize an objective function, such as generation cost, subject to both the network physics and engineering limits. The AC power flow equations introduce non-convexities in OPF problems. Due to these non-convexities, OPF problems may have multiple local optima\mrn{\cite{bukhsh_tps, NarimaniACC}} and are generally \mbox{NP-Hard}~\cite{bienstock2015nphard}. 

Many research efforts have focused on algorithms for obtaining locally optimal or approximate OPF solutions~\cite{ferc4}. Recent research has also developed convex relaxations of OPF problems~\cite{molzahn2017survey}. Convex relaxations bound the optimal objective values, can certify infeasibility, and, in some cases, provably provide globally optimal solutions to OPF problems.

The capabilities of convex relaxations are, in many ways, complementary to those of local solution algorithms. For instance, relaxations' objective value bounds can certify how close a local solution is to being globally optimal. Accordingly, local algorithms and relaxations are used together in spatial branch-and-bound methods~\cite{harsha2018pscc}. Solutions from relaxations are also useful for initializing some local solvers~\cite{marley2016}. Relaxations are also needed for certain solution algorithms for robust OPF problems~\cite{pscc2018robust}. Moreover, the objective value bounds provided by relaxations are directly useful in other contexts, e.g.,~\cite{molzahn_lesieutre_demarco-pfcondition,molzahn-opf_spaces}. The tractability and accuracy of these and other algorithms are largely determined by the employed relaxation's tightness. Tightening relaxations is thus an active research topic~\cite{molzahn2017survey}.


The quadratic convex (QC) relaxation is a promising approach that encloses the trigonometric and product terms in the polar representation of power flow equations within convex envelopes~\cite{coffrin2015qc}. These envelopes are formed with linear and second-order cone programming (SOCP) constraints, resulting in a convex formulation. The QC relaxation's tightness strongly depends on the quality of these convex envelopes. This paper focuses on improving these envelopes.

Previous work has proposed a variety of approaches for tightening the QC relaxation. These include valid inequalities, such as ``Lifted Nonlinear Cuts''~\cite{coffrin2016strengthen_tps,chen2015cuts} and constraints that exploit bounds on the differences in the voltage magnitudes~\cite{NarimaniPSCC2018}. Additionally, since the accuracies of the trigonometric and product envelopes in the QC relaxation rely on the voltage magnitude and angle difference bounds, bound tightening approaches can significantly strengthen the QC relaxation~\cite{coffrin2016strengthen_tps,chen2015,StrongSOCPRelaxations,arctan2,dmitry2019,Sundar}. When bound tightening approaches provide sign-definite angle difference bounds (i.e., the upper and lower bounds on the angle differences have the same sign), tighter trigonometric envelopes can be applied~\cite{coffrin2016strengthen_tps}.  

This paper proposes two improvements to further tighten QC relaxations of OPF problems. The first improvement leverages a polar representation of the branch admittances in addition to the rectangular representation used in previous QC formulations. 
%
%
Within certain ranges, portions of the trigonometric envelopes resulting from the polar admittance representation are at least as tight (and generally tighter) than the corresponding portions of the envelopes from the rectangular admittance representation. In other ranges, the trigonometric envelopes from the polar admittance representation neither contain nor are contained within the envelopes from the rectangular admittance representation. Thus, combining these envelopes tightens the QC relaxation, with empirical results suggesting limited impacts on solution times.


The polar admittance representation also enables our second improvement. We exploit a degree of freedom in the OPF formulation related to the per unit base power normalization. Selecting a \emph{complex base power} ($S_{base} = \left|S_{base}\right| e^{j\psi}$) results in a coordinate transformation that rotates the power flow equations relative to the typical choice of a real-valued base power. We leverage the associated rotational degree of freedom~$\psi$ to obtain tighter envelopes for the trigonometric functions. While previously proposed power flow algorithms~\cite{tortelli2015} and state estimation algorithms~\cite{ju2018} use similar formulations, this paper is, to the best of our knowledge, the first to exploit this rotational degree of freedom to improve convex relaxations.


This paper is organized as follows. Sections~\ref{OPF overview} and~\ref{Overview of QC Relaxation} review the OPF formulation and the previously proposed QC relaxation, respectively. Section~\ref{coordinate trasformation for OPF problem} describes the coordinate changes underlying our improved QC relaxation. Section~\ref{QC relaxation} then presents these improvements. Section~\ref{Numerical_results} empirically evaluates our approach. Section~\ref{conclusion} concludes the paper. 

\section{Overview of the Optimal Power Flow Problem}
\label{OPF overview}
This section formulates the OPF problem using a polar voltage phasor representation. The sets of buses, generators, and lines are $\mathcal{N}$, $\mathcal{G}$, and $\mathcal{L}$, respectively. \mrn{The set $\mathcal{R}$ contains the index of the bus that sets the angle reference.} Let $S_i^d = {P}_i^d + j{Q}_i^d$ and $S_i^g = {P}_i^g + j {Q}_i^g$ represent the complex load demand and generation, respectively, at bus~$i\in\mathcal{N}$, where $j = \sqrt{-1}$. Let $V_i$ and $\theta_i$ represent the voltage magnitude and angle at bus~$i\in\mathcal{N}$. Let $g_{sh,i} + j b_{sh,i}$ denote the shunt admittance at bus~$i\in\mathcal{N}$. For each generator, define a quadratic cost function with coefficients $c_{2,i} \geq 0$, $c_{1,i}$, and $c_{0,i}$.  For simplicity, we consider a single generator at each bus by setting the generation limits at buses without generators to zero. Upper and lower bounds for all variables are indicated by $\left(\overline{\,\cdot\,}\right)$ and $\left(\underline{\,\cdot\,}\right)$, respectively.

For ease of exposition, each line $\left(l,m\right)\in\mathcal{L}$ is modeled as a $\Pi$ circuit with mutual admittance $g_{lm}+j b_{lm}$ and shunt admittance $j b_{c,lm}$. Extensions to more general line models that allow for off-nominal tap ratios and non-zero phase shifts are straightforward and available in%
\ifarxiv 
~Appendix~\ref{appendix:power flow formulation with transformer}. 
\fi
Define $\theta_{lm}=\theta_{l}-\theta_{m}$ for $(l,m)\in\mathcal{L}$. The complex power flow into each line terminal $(l,m)\in\mathcal{L}$ is denoted by ${P}_{lm}+j{Q}_{lm}$, and the apparent power flow limit is $\overline{S}_{lm}$. 
The OPF problem is
\begin{subequations}
\label{OPF formulation}
\begin{align}
\label{eq:objective}
& \!\!\!\!\!\!\!\min\quad \sum_{{i}\in \mathcal{G}}
c_{2,i}\left(P_i^g\right)^2+c_{1,i}\,P_i^g+c_{0,i}\\ 
&\!\!\!\!\!\!\!\nonumber \text{subject to} \quad \left(\forall i\in\mathcal{N},\; \forall \left(l,m\right) \in\mathcal{L}\right) \\
\label{eq:pf1}
&\!\!\!\!\!\!\! P_i^g-P_i^d = g_{sh,i}\, V_i^2+\sum_{\substack{(l,m)\in \mathcal{L},\\\text{s.t.} \hspace{3pt} l=i}} \!P_{lm}+\!\!\sum_{\substack{(l,m)\in \mathcal{L},\\\text{s.t.} \hspace{3pt} m=i}} \!\!P_{ml}, \\
\label{eq:pf2}
&\!\!\!\!\!\!\! Q_i^g-Q_i^d = -b_{sh,i}\, V_i^2+\!\!\!\!\!\!\sum_{\substack{(l,m)\in \mathcal{L},\\ \text{s.t.} \hspace{3pt} l=i}} \!\!Q_{lm}+\!\!\!\sum_{\substack{(l,m)\in \mathcal{L},\\ \text{s.t.} \hspace{3pt} m=i}} \!\!\!\!Q_{ml},\\
\label{eq:angref}
&\!\!\!\!\!\!\! \mrn{\theta_r=0, \quad r\in\mathcal{R}}, \\
\label{eq:active_gen_limits}
&\!\!\!\!\!\!\! \underline{P}_i^g\leq P_i^g\leq \overline{P}_i^g, \quad 
\underline{Q}_i^g\leq Q_i^g\leq \overline{Q}_i^g,\\
\label{eq:bus_magnitude_limits}
&\!\!\!\!\!\!\! \underline{V}_i\leq V_{i} \leq \overline{V}_i,\\
\label{eq:bus_angle_limits}
&\!\!\!\!\!\!\!\underline{\theta}_{lm}\leq \theta_{lm}\leq \overline{\theta}_{lm},\\
\label{eq:pik}
&\!\!\!\!\!\!\! P_{lm} \!=\! g_{lm} V_l^2\! -\! g_{lm} V_l V_m\cos\left(\theta_{lm}\right)\! -\! b_{lm} V_l V_m\sin\left(\theta_{lm}\right),\\
\label{eq:qik}
&\!\!\!\!\!\!\! \nonumber Q_{lm} = -\left(b_{lm}+b_{c,lm}/2\right) V_l^2 + b_{lm} V_l V_m\cos\left(\theta_{lm}\right)\\ &\qquad\qquad  - g_{lm} V_l V_m\sin\left(\theta_{lm}\right),\\
\label{eq:pki}
&\!\!\!\!\!\!\!P_{ml}\! =\! g_{lm} V_m^2\! -\! g_{lm} V_l V_m\cos\left(\theta_{lm}\right)\! +\! b_{lm} V_l V_m\sin\left(\theta_{lm}\right),\\
\label{eq:qki}
&\!\!\!\!\!\!\!\nonumber Q_{ml} = -\left(b_{lm}+b_{c,lm}/2\right) V_m^2 + b_{lm} V_l V_m\cos\left(\theta_{lm}\right)\\ &\qquad\qquad  + g_{lm} V_l V_m\sin\left(\theta_{lm}\right),\\
\label{eq: apparent power limit sending}
&\!\!\!\!\!\!\! \left(P_{lm}\right)^2+\left(Q_{lm}\right)^2 \leq \left(\overline{S}_{lm}\right)^2, \quad 
\left(P_{ml}\right)^2+\left(Q_{ml}\right)^2 \leq \left(\overline{S}_{lm}\right)^2.
\end{align}
\end{subequations}
%
The objective~\eqref{eq:objective} minimizes the generation cost. Constraints~\eqref{eq:pf1} and~\eqref{eq:pf2} enforce power balance at each bus. Constraint~\eqref{eq:angref} sets the reference bus angle. The constraints in~\eqref{eq:active_gen_limits} bound the active and reactive power generation at each bus. Constraints~\eqref{eq:bus_magnitude_limits}--\eqref{eq:bus_angle_limits}, respectively, bound the voltage magnitudes and voltage angle differences. Constraints~\eqref{eq:pik}--\eqref{eq:qki} relate the active and reactive power flows with the voltage phasors at the terminal buses. The constraints in~\eqref{eq: apparent power limit sending} limit the apparent power flows into both terminals of each line.
\section{The QC Relaxation of the OPF Problem}
\label{Overview of QC Relaxation}
The QC relaxation convexifies the OPF problem~\eqref{OPF formulation} by enclosing the nonconvex terms in convex envelopes~\cite{coffrin2015qc}. The relevant nonconvex terms are the square $V_i^2$, $\forall i \in\mathcal{N}$, and the products $V_l V_m \cos(\theta_{lm})$ and $V_l V_m \sin(\theta_{lm})$, $\forall (l,m)\in\mathcal{L}$. The envelope for the generic squared function $x^2$ is $\langle x^2\rangle^T$:
\begin{align}
\label{eq:squareenvelopes}
\langle x^2\rangle^T =
\begin{cases}
\widecheck{x}: \begin{cases}\check{x} \geq x^2,\\
\widecheck{x} \leq \left({\overline{x}+\underline{x}}\right) x-{\overline{x}\underline{x}}.\\
\end{cases}
\end{cases}
\end{align}
where $\widecheck{x}$ is a lifted variable representing the squared term. 
Envelopes for the generic trigonometric functions $\sin(x)$ and $\cos(x)$ are $\left\langle \sin(x)\right\rangle^S$ and $\left\langle \cos(x)\right\rangle^C$:
\label{eq:convex_envelopes_sin&cos}
\begin{small}
\begin{align}
\label{eq:sine envelope}
 &\left\langle \sin(x)\right\rangle^S\!\!\! =\!
\begin{cases}
\!\widecheck{S}\!:\!\begin{cases}
\widecheck{S}\!\leq\!\cos\left(\frac{x^m}{2}\right)\!\!\left(x\!-\!\frac{x^m}{2}\right)\!\!+\!\sin \left(\frac{x^m}{2}\right) \text{if~} \underline{x}\!\leq0\!\le\overline{x},\\
\widecheck{S}\!\geq\!\cos\left(\frac{x^m}{2}\right)\!\!\left(x\!+\!\frac{x^m}{2}\right)\!\!-\!\sin\left(\frac{x^m}{2}\right)\text{if~} \underline{x}\!\leq0\!\le\overline{x},\\
\widecheck{S}\geq\frac{\sin\left({\underline{x}}\right)-\sin\left(\overline{x}\right)}{{\underline{x}-\overline{x}}}\left(x-{\underline{x}}\right)+\sin\left({\underline{x}}\right) \text{if~} \underline{x}\geq0,\\
\widecheck{S}\leq\frac{\sin\left({\underline{x}}\right)-\sin\left({\overline{x}}\right)}{{\underline{x}-\overline{x}}}\left(x-{\underline{x}}\right)+\sin\left({\underline{x}}\right) \text{if~} {\overline{x}}\leq0,
\end{cases}
\end{cases}\\[-0.3em]
\raisetag{18pt} \label{eq:cosine envelope}
&\left\langle\cos(x)\right\rangle^C =
\begin{cases}
\widecheck{C}:\!\begin{cases}
\widecheck{C}\leq 1-\frac{1-\cos\left({x^m}\right)}{\left(x^m\right)^2}x^2,\\
\widecheck{C}\geq\frac{\cos\left(\underline{x}\right)-\cos\left({\overline{x}}\right)}{{\underline{x}-\overline{x}}}\left(x-{\underline{x}}\right)+\cos\left({\underline{x}}\right),
\end{cases}
\end{cases}
\end{align}
\end{small}%
where $x^m = \max(\left|\underline{x}\right|,\left|\overline{x}\right|)$. \mrn{The envelopes $\left\langle \sin(x)\right\rangle^S$ and $\left\langle \cos(x)\right\rangle^C$ in~\eqref{eq:sine envelope} and~\eqref{eq:cosine envelope} are valid for $-\frac{\pi}{2} \leq x\leq \frac{\pi}{2}$. 
\ifarxiv
As detailed in Appendix~\ref{appendix: sine and cosine envelopes}, similar envelopes for the sine and cosine functions are defined analogously for other angle difference ranges.
\fi
} The lifted variables $\widecheck{S}$ and $\widecheck{C}$ are associated with the envelopes for the functions $\sin(\theta_{lm})$ and $\cos(\theta_{lm})$. 
%
%
%
%
The QC relaxation of the OPF problem in~\eqref{OPF formulation} is:
\begin{subequations}
\label{eq:QC relaxation}
\begin{align}
& \min \quad \sum_{\footnotesize{{i}\in \mathcal{N}}} c_{2,i}\left(P_i^g\right)^2+c_{1,i}\,P_i^g+c_{0,i}\\ 
&\nonumber \text{subject to} \quad \left(\forall i\in\mathcal{N},\; \forall \left(l,m\right) \in\mathcal{L}\right) \\
\label{eq:qc_p}
& P_i^g-P_i^d = g_{sh,i}\, w_{ii}+\!\!\!\!\!\sum_{\footnotesize{\substack{(l,m)\in \mathcal{L},\\ \text{s.t.} \hspace{3pt} l=i}}} \!\!\!P_{lm}+\sum_{\footnotesize{\substack{(l,m)\in \mathcal{L},\\ \text{s.t.} \hspace{3pt} m=i}}} \!\!P_{ml}, \\
\label{eq:qc_q}
& Q_i^g-Q_i^d = -b_{sh,i}\, w_{ii}+\!\!\!\sum_{\footnotesize{\substack{(l,m)\in \mathcal{L},\\ \text{s.t.} \hspace{3pt} l=i}}} \!\!Q_{lm}+\sum_{\footnotesize{\substack{(l,m)\in \mathcal{L},\\ \text{s.t.} \hspace{3pt} m=i}}} \!\!Q_{ml},\\
\label{eq:angref_QC}
& \mrn{\theta_r=0, \quad r\in\mathcal{R}}, \\
\label{eq:active_gen_limits_QC}
& \mrn{\underline{P}_i^g\leq P_i^g\leq \overline{P}_i^g, \quad 
\underline{Q}_i^g\leq Q_i^g\leq \overline{Q}_i^g,}\\
\label{eq:bus_magnitude_limits_QC}
& \mrn{\underline{V}_i\leq V_{i} \leq \overline{V}_i},\\
\label{eq:bus_angle_limits_QC}
&\mrn{\underline{\theta}_{lm}\leq \theta_{lm}\leq \overline{\theta}_{lm}},\\
\label{eq:qc_V}
&  (\underline{V}_i)^2\leq w_{ii} \leq (\overline{V}_i)^2, \qquad w_{ii} \in\left\langle V_i^2 \right\rangle^T,\\
\label{eq:qc_pik}
& P_{lm} = g_{lm} w_{ll} - g_{lm} c_{lm} - b_{lm} s_{lm},\\
\label{eq:qc_qik}
& Q_{lm} = -\left(b_{lm}+b_{c,lm}/2\right) w_{ll} + b_{lm} c_{lm}- g_{lm} s_{lm}, \\
\label{eq:qc_pki}
& P_{ml} = g_{lm} w_{mm} - g_{lm} c_{lm} + b_{lm} s_{lm}, \\
\label{eq:qc_qki}
& Q_{ml} = -\left(b_{lm}+b_{c,lm}/2\right) w_{mm} + b_{lm} c_{lm}+ g_{lm} s_{lm},\\
\label{eq: apparent power limit sending_QC}
&\mrn{ \left(P_{lm}\right)^2+\left(Q_{lm}\right)^2 \leq \left(\overline{S}_{lm}\right)^2, \quad 
\left(P_{ml}\right)^2+\left(Q_{ml}\right)^2 \leq \left(\overline{S}_{lm}\right)^2,}\\
    \label{eq:qc_cik}
    \nonumber & {c}_{lm}  = \sum_{\footnotesize{k=1,\ldots,8}}\!\! \lambda_k\, \rho^{(k)}_1 \rho^{(k)}_2 \rho^{(k)}_3,  \quad \widecheck{C}_{lm} \in \left\langle\cos(\theta_{lm})\right\rangle^C,\\
    \nonumber & V_l =\!\!\!\! \sum_{\footnotesize{k=1,\ldots,8}}\!\! \lambda_k\rho^{(k)}_1,\; V_m =\!\!\!\! \sum_{\footnotesize{k=1,\ldots,8}}\!\! \lambda_k\rho^{(k)}_2,\; \widecheck{C}_{lm} =\!\!\!\! \sum_{\footnotesize{k=1,\ldots,8}}\!\! \lambda_k\rho^{(k)}_3,\\
     & \sum_{k=1,\ldots,8} \lambda_k = 1, \quad \lambda_k \geqslant 0, \quad k=1,\ldots,8.\\
        \label{eq:qc_sik}
          \nonumber & {s}_{lm}  = \sum_{\footnotesize{k=1,\ldots,8}}\!\!\! \gamma_k\, \zeta^{(k)}_1 \zeta^{(k)}_2 \zeta^{(k)}_3, \quad \widecheck{S}_{lm} \in \left\langle\sin(\theta_{lm})\right\rangle^S,\\
    \nonumber & V_l =\!\!\!\! \sum_{\footnotesize{k=1,\ldots,8}}\!\! \gamma_k\zeta^{(k)}_1,\; V_m =\!\!\!\! \sum_{\footnotesize{k=1,\ldots,8}}\!\! \gamma_k\zeta^{(k)}_2,\; \widecheck{S}_{lm} =\!\!\!\! \sum_{\footnotesize{k=1,\ldots,8}}\!\!\ \gamma_k\zeta^{(k)}_3,\\
     & \sum_{\footnotesize{k=1,\ldots,8}} \gamma_k = 1, \quad \gamma_k \geqslant 0, \quad k=1,\ldots,8.  \\
\label{eq:ell_apparent}
& P_{lm}^2 + Q_{lm}^2 \leq \mrn{w_{ll}}\, \ell_{lm},\\
\label{eq:ell_expression}
&\ell_{lm} = \left(Y_{lm}^2- \frac{b_{c,lm}^2}{4}\right)\mrn{w_{ll}} + Y_{lm}^2 \mrn{w_{mm}} - 2Y_{lm}^2c_{lm} - b_{c,lm}Q_{lm},\\
\label{eq:qc_linking}
& \mrn{\begin{bmatrix}
\lambda_1+\lambda_2-\gamma_1-\gamma_2 \\
\lambda_3+\lambda_4-\gamma_3-\gamma_4 \\
\lambda_5+\lambda_6-\gamma_5-\gamma_6 \\
\lambda_7+\lambda_8-\gamma_7-\gamma_8
\end{bmatrix}^\intercal
\begin{bmatrix}
\underline{V}_l \underline{V}_m\\
\underline{V}_l \overline{V}_m \\
\overline{V}_l \underline{V}_m \\
\overline{V}_l \overline{V}_m
\end{bmatrix}=0}.
\end{align}
\end{subequations}
where $\ell_{lm}$ represents the squared magnitude of the current flow into terminal $l$ of line $(l,m)\in\mathcal{L}$ \mrn{and $\left(\cdot\right)^\intercal$ is the transpose operator}. The relationship between $\ell_{lm}$ and the power flows $P_{lm}$ and $Q_{lm}$ in~\eqref{eq:ell_apparent} tightens the QC relaxation~\cite{farivar2011,coffrin2015qc}.
\ifarxiv
~Appendix~\ref{appendix:power flow formulation with transformer} gives an expression for $\ell_{lm}$ that considers lines with off-nominal tap ratios and non-zero phase shifts. 
\fi
Also, as shown in~\eqref{eq:qc_V}, $w_{ii}$ is associated with the squared voltage magnitude at bus~$i$.

The lifted variables $c_{lm}$ and $s_{lm}$ represent relaxations of the trilinear terms $V_lV_m\cos(\theta_{lm})$ and $V_lV_m\sin(\theta_{lm})$, respectively, with~\eqref{eq:qc_cik} and~\eqref{eq:qc_sik} formulating an ``extreme point'' representation of the convex hulls for the trilinear terms $V_l V_m \widecheck{C}_{lm}$ and $V_l V_m \widecheck{S}_{lm}$\mrn{\cite{harsha2018pscc, GlobalSIP2018}}. 
The auxiliary variables $\lambda_k,\; \gamma_k\in[0,1]$, $k=1,\ldots,8$, are used in the formulations of these convex hulls.
The extreme points of $V_lV_m\widecheck{C}_{lm}$ are $\rho^{(k)} \in [\underline{V_l},\overline{V_l}] \times [\underline{V_m},\overline{V_m}]\times[\underline{\widecheck{C}}_{lm},\overline{\widecheck{C}}_{lm}]$, $k=1,\ldots,8$ and the extreme points of $V_l V_m\widecheck{S}_{lm}$ are $\zeta^{(k)} \in [\underline{V_l},\overline{V_l}] \times [\underline{V_m},\overline{V_m}]\times[\underline{\widecheck{S}}_{lm},\overline{\widecheck{S}}_{lm}]$, $k=1,\ldots,8$. Since sine and cosine are odd and even functions, respectively, $c_{lm} = c_{ml}$ and $s_{lm} = -s_{ml}$.

A ``linking constraint'' from~\cite{Sundar} is also enforced in~\eqref{eq:qc_linking}. This linking constraint is associated with the bilinear terms $V_lV_m$ that are shared in $V_lV_m\cos(\theta_{lm})$ and $V_lV_m\sin(\theta_{lm})$.

\section{Coordinate Transformations}
\label{coordinate trasformation for OPF problem}
The improvements to the QC relaxation's envelopes that are our main contributions are based on certain coordinate transformations. This section describes these transformations. We first form the power flow equations using polar representations of the lines' mutual admittances. We then introduce a complex base power in the per unit normalization that provides a rotational degree of freedom in the power flow equations.

While this section uses a $\Pi$ circuit line model for simplicity, extensions to more general line models are straightforward. These extensions are presented in%
\ifarxiv
~Appendix~\ref{appendix:power flow formulation with transformer}.
\fi

\subsection{Power Flow Equations with Admittance in Polar Form}
\label{Power flow with Admittance in Polar Form}
Equations~\eqref{eq:pik}--\eqref{eq:qki} model the power flows through a line~$(l,m)\in\mathcal{L}$ via a rectangular representation of the line's mutual admittance, $g_{lm} + jb_{lm}$. In~\eqref{eq:qc_pik}--\eqref{eq:qc_qki}, the QC relaxation from~\cite{coffrin2015qc} uses this rectangular admittance representation.

The line flows can be equivalently modeled using a polar representation of the mutual admittance, $Y_{lm}e^{j\delta_{lm}}$, where $Y_{lm} = \sqrt{g_{lm}^2 + b_{lm}^2}$ and $\delta_{lm} = \arctan{\left(b_{lm}/g_{lm}\right)}$ are the magnitude and angle of the mutual admittance for line $(l,m)\in\mathcal{L}$, respectively. Using polar admittance coordinates, the complex power flows $S_{lm}$ and $S_{ml}$ into each line terminal are:
\begin{subequations}
\small
\label{complex power flow}
\begin{align}
\label{complex power flow sending}
 &S_{lm} = V_l e^{j\theta_l}
  \left(\left(Y_{lm}e^{j\delta_{lm}} + j\frac{b_{c,lm}}{2} \right)V_l e^{j\theta_l}  - Y_{lm} e^{j\delta_{lm}} V_m e^{j\theta_m} \right)^\ast\!\!\!\!,\\
  \label{complex power flow receiving}
 &S_{ml} = V_m e^{j\theta_m}
  \!\left(\!-Y_{lm}e^{j\delta_{lm}}V_l e^{j\theta_l}\!\!+ \!\!\left(Y_{lm}e^{j\delta_{lm}}\!+ \!j\frac{b_{c,lm}}{2}\!\right) V_m e^{j\theta_m}\!\!\right)^\ast\!\!\!\!,
\end{align}
\end{subequations}
where $(\cdot)^\ast$ is the complex conjugate.
Taking the real and imaginary parts of~\eqref{complex power flow} yields the active and reactive line flows: 
\begin{subequations}
\label{power flows with polar representation of admittance}
{\small
\begin{align} 
   \label{eq: polar pik} & P_{lm} = \Re(S_{lm}) = Y_{lm}\cos(\delta_{lm})V_l^2 -Y_{lm}V_lV_m \cos(\theta_{lm}-\delta_{lm}),\\[-0.8em]
   \nonumber & Q_{lm} = \Im(S_{lm}) =  -\left(Y_{lm}\sin(\delta_{lm})+b_{c,lm}/2\right)V_l^2 \\ \label{eq: polar qik} & \qquad\qquad\qquad\qquad\qquad -  Y_{lm}V_lV_m \sin(\theta_{lm}-\delta_{lm}),\\
   \label{eq: polar pki}
    & P_{ml} = \Re(S_{ml}) = Y_{lm}\cos(\delta_{lm})V_m^2 - Y_{lm}V_lV_m \cos(\theta_{lm}+\delta_{lm}),\\[-0.8em]
   \nonumber & Q_{ml} = \Im(S_{ml}) = -\left(Y_{lm}\sin(\delta_{lm})+b_{c,lm}/2\right)V_m^2 \\ 
   \label{eq: polar qki} & \qquad\qquad\qquad\qquad\qquad + Y_{lm}V_lV_m\sin(\theta_{lm}+\delta_{lm}).
\end{align}}%
\end{subequations}

With the rectangular admittance representation, the active and reactive power flow equations~\eqref{eq:pik}--\eqref{eq:qik} each have two trigonometric terms (i.e., $\cos(\theta_{lm})$ and $\sin(\theta_{lm})$). Conversely, there is only one trigonometric term in each of the power flow equations that use the polar admittance representation~\eqref{power flows with polar representation of admittance} (e.g., $\cos(\theta_{lm}-\delta_{lm})$ for $P_{lm}$ and $\sin(\theta_{lm}-\delta_{lm})$ for $Q_{lm}$). While these formulations are equivalent, the differing representations of the trigonometric terms suggest the possibility of using different trigonometric envelopes. The QC formulation we will propose in Section~\ref{QC relaxation for Rotated base OPF} exploits these differences.
\subsection{Rotated Power Flow Formulation}
\label{Rotated Base OPF formulation}
The base power used in the per unit normalization is traditionally chosen to be a real-valued quantity. More generally, complex-valued choices for the base power are also acceptable and can provide benefits for some algorithms. For instance, certain power flow~\cite{tortelli2015} and state estimation algorithms~\cite{wei2015,ju2018} leverage formulations with a complex-valued base power. 

To improve the QC relaxation's trigonometric envelopes, this section reformulates the OPF problem with a complex base power. Let $S_{base}^{orig}$ and $S_{base}^{new} e^{j\psi}$ denote the original and the new base power, respectively, where $S_{base}^{orig}$, $S_{base}^{new}$, and $\psi$ are real-valued. Thus, the original base $S_{base}^{orig}$ is real-valued, while the new base $S_{base}^{new}e^{j\psi}$ is complex-valued with magnitude $S_{base}^{new}$ and angle $\psi$. Quantities associated with the new base power will be accented with a tilde, $\left(\tilde{\,\cdot\,}\right)$. Complex power flows in the original base and the new base are related as:
\begin{equation*}
\tilde{S}_{lm} = S_{lm}\cdot \frac{S_{base}^{orig}}{S_{base}^{new}e^{j\psi}}, \qquad \tilde{S}_{ml} = S_{ml}\cdot \frac{S_{base}^{orig}}{S_{base}^{new}e^{j\psi}}.
\end{equation*}
Since changing the magnitude of the base power does not affect the arguments of the trigonometric functions in the power flow equations, we choose $S_{base}^{new} = S_{base}^{old}$. With this~choice,
\begin{equation*}
\tilde{S}_{lm} = S_{lm}/e^{j\psi}, \qquad \tilde{S}_{ml} = S_{ml}/e^{j\psi}.
\end{equation*}
The angle of the base power, $\psi$, affects the arguments of the trigonometric functions, as shown in the following derivation:
\begin{subequations}
\vspace*{-1em}
\label{rotated power flow complex}
\begin{align}
\label{rotated power flow sending}
&\nonumber \tilde{S}_{lm}=S_{lm}/e^{j\psi}= \left(Y_{lm}e^{-j(\delta_{lm}+\psi)}+(b_{c,lm}/2)e^{-j(\frac{\pi}{2}+\psi)}\right)V_l^2\\&\qquad\qquad\qquad\qquad - Y_{lm}V_l V_m e^{j(-\delta_{lm}+\theta_{lm}-\psi)},\\
  \label{rotated power flow receiving}
&\nonumber \tilde{S}_{ml}=S_{ml}/e^{j\psi} = \left(Y_{lm}  e^{-j(\delta_{lm}+\psi)}+(b_{c,lm}/2)  e^{-j(\frac{\pi}{2}+\psi)}\right)V_m^2\\&\qquad\qquad\qquad\qquad - Y_{lm}V_m V_l e^{-j(\delta_{lm}+ \theta_{lm}+\psi)}.
\end{align}
\end{subequations}
Taking the real and imaginary parts of~\eqref{rotated power flow complex} yields:
\begin{subequations}
\interdisplaylinepenalty=10000
\label{rotated power flow}
\begin{align}
\nonumber & \tilde{P}_{lm}\! =\! \Re(\tilde{S}_{lm})\! =\! \left(Y_{lm}\cos(\delta_{lm}+\psi)-(b_{c,lm}/2)\sin(\psi)\right)V_l^2\\ \label{eq:shifted OPF active power flow lm} &\qquad\quad\quad\qquad\qquad -
  Y_{lm}V_l V_m\cos(\theta_{lm}-\delta_{lm}-\psi),\\
  \label{eq:shifted OPF reactive power flow lm}
\nonumber & \tilde{Q}_{lm} \!=\! \Im(\tilde{S}_{lm}) \!=\! -\left(Y_{lm}\sin(\delta_{lm}+\psi)\!+\!(b_{c,lm}/2)\cos(\psi)\right)\!V_l^2\\ &\qquad\quad\quad\qquad\qquad -
  Y_{lm}V_l V_m\sin(\theta_{lm}-\delta_{lm}-\psi),\\
  \label{eq:shifted OPF active power flow ml}
\nonumber & \tilde{P}_{ml} \!=\! \Re(\tilde{S}_{ml}) \!=\! \left(Y_{lm} \cos(\delta_{lm}+\psi)\!-\!(b_{c,lm}/2)  \sin(\psi)\right)\!V_m^2
  \\ &\qquad\quad\quad\qquad\qquad - Y_{lm}V_m V_l \cos(\theta_{lm}+\delta_{lm}+\psi),\\
\nonumber & \tilde{Q}_{ml} \!=\! \Im(\tilde{S}_{ml}) \!=\! - \left(Y_{lm}  \sin(\delta_{lm}+\psi)\!+\!(b_{c,lm}/2) \cos(\psi)\right)\!V_m^2 \\ &\qquad\quad\quad\qquad\qquad + Y_{lm}V_m V_l \sin(\theta_{lm}+\delta_{lm}+\psi).
  \label{eq:shifted OPF reactive power flow ml}
\end{align}
\end{subequations}
The arguments of the trigonometric functions $\cos(\theta_{lm}-\delta_{lm}-\psi)$, $\sin(\theta_{lm}-\delta_{lm}-\psi)$, $\cos(\theta_{lm}+\delta_{lm}+\psi)$, and $\sin(\theta_{lm}+\delta_{lm}+\psi)$ in~\eqref{rotated power flow} are linear in $\psi$. For a given $\psi$, all other trigonometric terms in~\eqref{rotated power flow} are constants that do not require special handling.
\vspace{-.4cm}
\subsection{Rotated OPF Problem}

We next represent the complex power generation and load demands using the new base power:
\begin{equation*}
\tilde{S}_i^g = S_i^g\cdot \frac{S_{base}^{orig}}{S_{base}^{new}e^{j\psi}} = \frac{S_i^g}{e^{j\psi}} = \frac{P_i^g + j Q_i^g}{e^{j\psi}}.
\end{equation*}
Define $\tilde{S}_i^g = \tilde{P}_i^g + j\tilde{Q}_i^g$, $\forall i\in\mathcal{N}$. Taking the real and imaginary parts of $\tilde{S}_i^g$ yields the following relationship between the power generation in the new and original bases:
%
\begin{align}
\label{eq:gen transformation}
& \begin{bmatrix}
\tilde{P}_i^g \\
\tilde{Q}_i^g
\end{bmatrix} = \begin{bmatrix}
\cos(\psi) & \sin(\psi) \\
-\sin(\psi) & \cos(\psi)
\end{bmatrix}
\begin{bmatrix}
P_i^g \\
Q_i^g
\end{bmatrix}.
\end{align}
The inverse relationship is well defined for any choice of $\psi$ since the matrix in~\eqref{eq:gen transformation} is invertible. 

The analogous relationship for the power demands is:
\begin{align}
\label{eq:load transformation}
& \begin{bmatrix}
\tilde{P}_i^d \\
\tilde{Q}_i^d
\end{bmatrix} = \begin{bmatrix}
\cos(\psi) & \sin(\psi) \\
-\sin(\psi) & \cos(\psi)
\end{bmatrix}
\begin{bmatrix}
P_i^d \\
Q_i^d
\end{bmatrix}.
\end{align}

Applying~\eqref{rotated power flow}--\eqref{eq:load transformation} to~\eqref{OPF formulation} yields a ``rotated'' OPF problem:
\begin{subequations}
\ifarxiv
\interdisplaylinepenalty=10000
\fi
\label{shifted OPF}
\begin{align}
\label{eq:shifted OPF objective function}
&\nonumber \min\quad \textstyle\sum_{\footnotesize{{i}\in \mathcal{G}}}
c_{2,i}\left( {\tilde P_i^g\cos (\psi ) - \tilde Q_i^g\sin (\psi)} \right)^2\\[-0.5em]
&\quad\qquad\qquad+c_{1,i}\left( {\tilde P_i^g\cos (\psi ) - \tilde Q_i^g\sin (\psi )} \right)+c_{0,i}\\
&\nonumber \text{subject to} \quad \left(\forall i\in\mathcal{N},\; \forall \left(l,m\right) \in\mathcal{L}\right)\\ 
\nonumber
& \tilde{P}_i^g-\tilde{P}_i^d =\left(g_{sh,i}  \cos(\psi)-b_{sh,i} \sin(\psi)\right)V_i^2\\ \label{eq:shifted OPF active power injection} & \quad\quad\qquad\qquad + \sum_{\substack{(l,m)\in \mathcal{L},\\ \text{s.t.} \hspace{3pt} l=i}} P_{lm}+\sum_{\substack{(l,m)\in \mathcal{L},\\ \text{s.t.} \hspace{3pt} m=i}} P_{ml}, \\
\label{eq:shifted OPF reactive power injection}
& \tilde{Q}_i^g-\tilde{Q}_i^d =-\left(g_{sh,i}\sin(\psi)+b_{sh,i}  \cos(\psi)\right)V_i^2\\ \nonumber&\quad\quad\qquad\qquad +\sum_{\substack{(l,m)\in \mathcal{L},\\ \text{s.t.} \hspace{3pt} l=i}} Q_{lm}+\sum_{\substack{(l,m)\in \mathcal{L},\\ \text{s.t.} \hspace{3pt} m=i}} Q_{ml},\\[-0.5em]
\label{eq:shifted OPF slack bus}
& \mrn{\theta_r=0, \quad r\in\mathcal{R}},\\
\label{eq:shifted OPF active power limits}
& \underline{P}_i^{g}\leq \tilde P_i^g\cos (\psi ) - \tilde Q_i^g\sin (\psi )\leq \overline{P}_i^{g},\\
 \label{eq:shifted OPF reactive power limits}
& \underline{Q}_i^{g}\leq \tilde Q_i^g\cos (\psi ) + \tilde P_i^g\sin (\psi )\leq \overline{Q}_i^{g},\\
 \label{eq:shifted OPF voltage limits}
& \underline{V}_i\leq V_{i} \leq \overline{V}_i, \quad\qquad
\underline{\theta}_{lm}\leq \theta_{lm}\leq \overline{\theta}_{lm},\\
  \label{eq:shifted OPF power flow limit lm}
& (\tilde{P}_{lm})^2+(\tilde{Q}_{lm})^2 \leq (\overline{S}_{lm})^2,~~
 (\tilde{P}_{ml})^2+(\tilde{Q}_{ml})^2 \leq (\overline{S}_{lm})^2,\\
& \text{Eq.}~\eqref{rotated power flow}.
\end{align}
\end{subequations}
The rotated OPF problem~\eqref{shifted OPF} is equivalent to~\eqref{OPF formulation} in that any solution $\{ V^\ast,\theta^\ast,\tilde{P}^{g\star},\tilde{Q}^{g\star}\}$ to~\eqref{shifted OPF} can be mapped to a solution $\{ V^\ast,\theta^\ast,P^{g\star},Q^{g\star}\}$ to~\eqref{OPF formulation} using~\eqref{eq:gen transformation}. Solutions to both formulations have the same voltage magnitudes and angles, $V^\ast$ and $\theta^\ast$. Thus,~\eqref{shifted OPF} can be interpreted as revealing a degree of freedom associated with choosing the base power's phase angle $\psi$. The next section exploits this degree of freedom to tighten the QC relaxation's trigonometric envelopes.

\section{Rotated QC relaxation}
\label{QC relaxation}
This section leverages the coordinate transformations presented in Section~\ref{coordinate trasformation for OPF problem} to tighten the QC relaxation. We first propose and analyze new envelopes for the trigonometric functions and trilinear terms. We then describe an empirical analysis that informs the choice of the base power angle~$\psi$ in order to tighten the relaxation for typical OPF problems.

 
\subsection{Convex Envelopes for the Trigonometric Terms}
\label{Different envelope epresentation explanations}

A key determinant of the QC relaxation's tightness is the quality of the convex envelopes for the trigonometric terms in the power flow equations. The rotated OPF formulation~\eqref{shifted OPF} has four relevant trigonometric terms for each line: $\cos(\theta_{lm}-\delta_{lm}-\psi)$, $\sin(\theta_{lm}-\delta_{lm}-\psi)$, $\cos(\theta_{lm}+\delta_{lm}+\psi)$, and $\sin(\theta_{lm}+\delta_{lm}+\psi)$, $\forall (l,m)\in\mathcal{L}$. This contrasts with the two unique trigonometric terms ($\cos(\theta_{lm})$ and $\sin(\theta_{lm})$) per pair of connected buses in the OPF formulation~\eqref{OPF formulation}.

\mrn{This would seem to suggest that at least twice as many lifted variables would be required in order to relax the rotated OPF formulation~\eqref{shifted OPF} compared to the original OPF formulation~\eqref{OPF formulation}, i.e., $\tilde{C}_{lm}^{(s)}$, $\tilde{S}_{lm}^{(s)}$ and $\tilde{C}_{lm}^{(r)}$, $\tilde{S}_{lm}^{(r)}$ for relaxing the sending end quantities $\cos(\theta_{lm}-\delta_{lm}-\psi)$, $\sin(\theta_{lm}-\delta_{lm}-\psi)$ and the receiving end quantities $\cos(\theta_{lm}+\delta_{lm}+\psi)$, $\sin(\theta_{lm}+\delta_{lm}+\psi)$ in~\eqref{shifted OPF} versus $\widecheck{C}_{lm}$ and $\widecheck{S}_{lm}$ for relaxing $\cos(\theta_{lm})$ and $\sin(\theta_{lm})$ in~\eqref{OPF formulation}. However, this is not the case since the arguments of the trigonometric terms in the rotated OPF formulation are not independent.} For notational convenience, define $\hat{\delta}_{lm} = \delta_{lm} + \psi$. The angle sum and difference identities imply the following relationships:
\begin{equation}
\small{
\label{eq:angle_sumdiff}
\!\!\begin{bmatrix}
\sin(\hat{\delta}_{lm} + \theta_{lm}) \\
\cos(\hat{\delta}_{lm} + \theta_{lm}) \\
\sin(\hat{\delta}_{lm} - \theta_{lm}) \\
\cos(\hat{\delta}_{lm} - \theta_{lm})
\end{bmatrix} \!=\! \begin{bmatrix}
\sin(\hat{\delta}_{lm}) & \cos(\hat{\delta}_{lm}) \\
\cos(\hat{\delta}_{lm}) & -\sin(\hat{\delta}_{lm}) \\
\sin(\hat{\delta}_{lm}) & -\cos(\hat{\delta}_{lm}) \\
\cos(\hat{\delta}_{lm}) & \sin(\hat{\delta}_{lm})
\end{bmatrix}\begin{bmatrix}
\cos(\theta_{lm}) \\
\sin(\theta_{lm})
\end{bmatrix}.
}
\end{equation}
Rearranging these relationships yields:
\begin{equation}
\label{eq:terminal_relationships}
\begin{bmatrix}
\sin(\theta_{lm} + \hat{\delta}_{lm}) \\
\cos(\theta_{lm} + \hat{\delta}_{lm})
\end{bmatrix}\! = \!\begin{bmatrix}
\alpha_{lm} & \beta_{lm} \\
-\beta_{lm} & \alpha_{lm}
\end{bmatrix}
\begin{bmatrix}
\sin(\theta_{lm} - \hat{\delta}_{lm})\\
\cos(\theta_{lm} - \hat{\delta}_{lm})
\end{bmatrix}.
\end{equation}
where, for notational convenience, $\alpha_{lm} =(\cos(\hat{\delta}_{lm}))^2 - (\sin(\hat{\delta}_{lm}))^2$ and $\beta_{lm} = 2\cos(\hat{\delta}_{lm})\sin(\hat{\delta}_{lm})$. \mrn{The implication of~\eqref{eq:terminal_relationships} is that only two (rather than four) lifted variables are defined per line (chosen to be the sending end quantities $\tilde{C}_{lm}^{(s)}$ and $\tilde{S}_{lm}^{(s)}$ for relaxing the trigonometric terms $\cos(\theta_{lm} - \hat{\delta}_{lm})$ and $\sin(\theta_{lm} - \hat{\delta}_{lm})$).} The remaining trigonometric functions, $\sin(\theta_{lm} + \hat{\delta}_{lm})$ and $\cos(\theta_{lm} + \hat{\delta}_{lm})$, are representable in terms of $\sin(\theta_{lm} - \hat{\delta}_{lm})$ and $\cos(\theta_{lm} - \hat{\delta}_{lm})$ via the linear transformation~\eqref{eq:terminal_relationships}. Since the matrix in~\eqref{eq:terminal_relationships} is invertible for all $\hat{\delta}_{lm}$, the transformation~\eqref{eq:terminal_relationships} is always well-defined. 

\mrn{While not explicitly including the lifted variables $\tilde{C}_{lm}^{(r)}$ and $\tilde{S}_{lm}^{(r)}$ for the receiving end quantities, we tighten the relaxation of~\eqref{shifted OPF} by enforcing the trigonometric envelopes associated with both the sending and receiving end quantities using~\eqref{eq:terminal_relationships}:
\begin{subequations}
\label{eq:rotated_envelopes_sending}
\begin{align}
\tilde{C}_{lm}^{(s)} \in \left\langle\cos(\theta_{lm}-\delta_{lm}-\psi)\right\rangle^{C}, \\
\tilde{S}_{lm}^{(s)} \in \left\langle\sin(\theta_{lm}-\delta_{lm}-\psi)\right\rangle^{S},
\end{align}
\end{subequations}\vspace*{-2em}
\begin{subequations}
\label{eq:rotated_envelopes_receiving}
\begin{align}
\alpha_{lm} \tilde{S}_{lm}^{(s)} + \beta_{lm} \tilde{C}_{lm}^{(s)} \in \left\langle\sin(\theta_{lm}+\delta_{lm}+\psi)\right\rangle^{S}, \\
-\beta_{lm} \tilde{S}_{lm}^{(s)} + \alpha_{lm} \tilde{C}_{lm}^{(s)} \in \left\langle\cos(\theta_{lm}+\delta_{lm}+\psi)\right\rangle^{C},
\end{align}
\end{subequations}
}



Related special consideration is needed for parallel lines. While the rest of this section considers systems without parallel lines for simplicity,%
\ifarxiv
~Appendix~\ref{appendix:parallel}
\fi discusses this issue in detail. Using the linear relationships in~\eqref{eq:terminal_relationships}%
\ifarxiv ~(and in~\eqref{eq:terminal_relationships_parallel} from Appendix~\ref{appendix:parallel} for systems with parallel lines),
\fi all relevant trigonometric terms in~\eqref{shifted OPF} can be represented as linear combinations of $\sin(\theta_{lm} - \delta_{lm} - \psi)$ and $\cos(\theta_{lm} - \delta_{lm} - \psi)$ for each unique pair of connected buses $(l,m)\in\mathcal{L}$. 
\mrn{The QC relaxations of~\eqref{OPF formulation} and~\eqref{shifted OPF} hence have the same number of lifted variables (two per pair of connected buses).}

There are two characteristics that distinguish the trigonometric expressions in~\eqref{OPF formulation} and~\eqref{shifted OPF}: \mrn{First, the power flow equations~\eqref{eq:pik}--\eqref{eq:qki} contain weighted sums of two trigonometric functions of $\theta_{lm}$, while~\eqref{eq:shifted OPF active power flow lm}--\eqref{eq:shifted OPF reactive power flow ml} each contain a single trigonometric function of $\theta_{lm}$. (The trigonometric expressions $\cos(\delta_{lm} + \psi)$, $\sin(\delta_{lm} + \psi)$, $\cos(\psi)$, and $\sin(\psi)$ in~\eqref{eq:shifted OPF active power flow lm}--\eqref{eq:shifted OPF reactive power flow ml} are constants that do not require special consideration.)} Second, the base power angle $\psi$ used to formulate~\eqref{shifted OPF} provides a degree of freedom that shifts the arguments of the trigonometric functions in~\eqref{eq:shifted OPF active power flow lm}--\eqref{eq:shifted OPF reactive power flow ml}. We next discuss how both of these characteristics can be exploited to tighten the QC relaxation.

Regarding the first distinguishing characteristic, factoring out~$-V_l V_m$ to focus on the trigonometric functions shows that the relaxation of~\eqref{eq:pik} depends on the quality of a weighted sum of trigonometric envelopes: $g_{lm}\left\langle\cos\left(\theta_{lm}\right)\right\rangle^C + b_{lm}\left\langle\sin\left(\theta_{lm}\right)\right\rangle^S$. The relaxation of~\eqref{eq:shifted OPF active power flow lm} depends on the quality of the envelope $Y_{lm}\left\langle\cos\left(\theta_{lm}-\delta_{lm}-\psi \right)\right\rangle^C$. (The relaxations of~\eqref{eq:pik}--\eqref{eq:qki} and~\eqref{eq:shifted OPF active power flow lm}--\eqref{eq:shifted OPF reactive power flow ml} are analogous.) To focus on the first characteristic, consider the latter envelope with $\psi = 0$. 

\begin{figure}[p!]
\centering
\subfloat[$-90^\circ \leq \theta_{lm} \leq 90^\circ$]{\includegraphics[width=6.5cm]{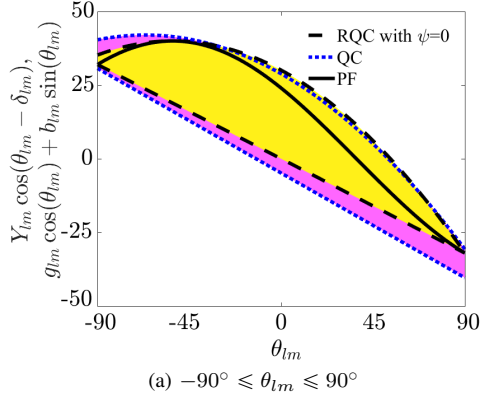}}\\ 
   \subfloat[$-60^\circ \leq \theta_{lm} \leq 0^\circ$]{\includegraphics[width=6.5cm]{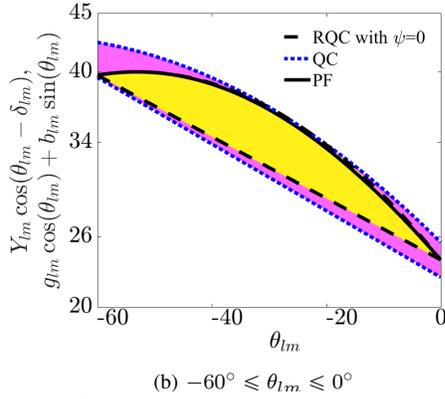}}\\   \vspace{-.3cm}
   \vspace{-.1cm}
   \subfloat[$-30^\circ \leq \theta_{lm} \leq 30^\circ$]{\includegraphics[width=6.5cm]{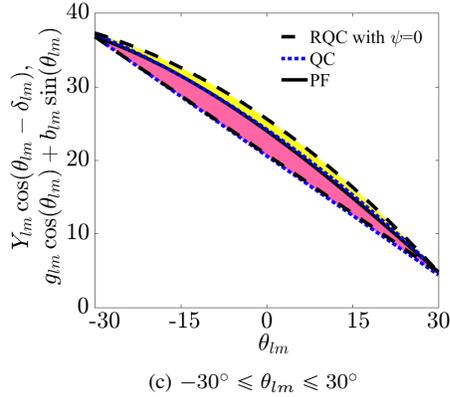}}\\
   \subfloat[$-30^\circ \leq \theta_{lm} \leq 0^\circ$]{\includegraphics[width=6.5cm]{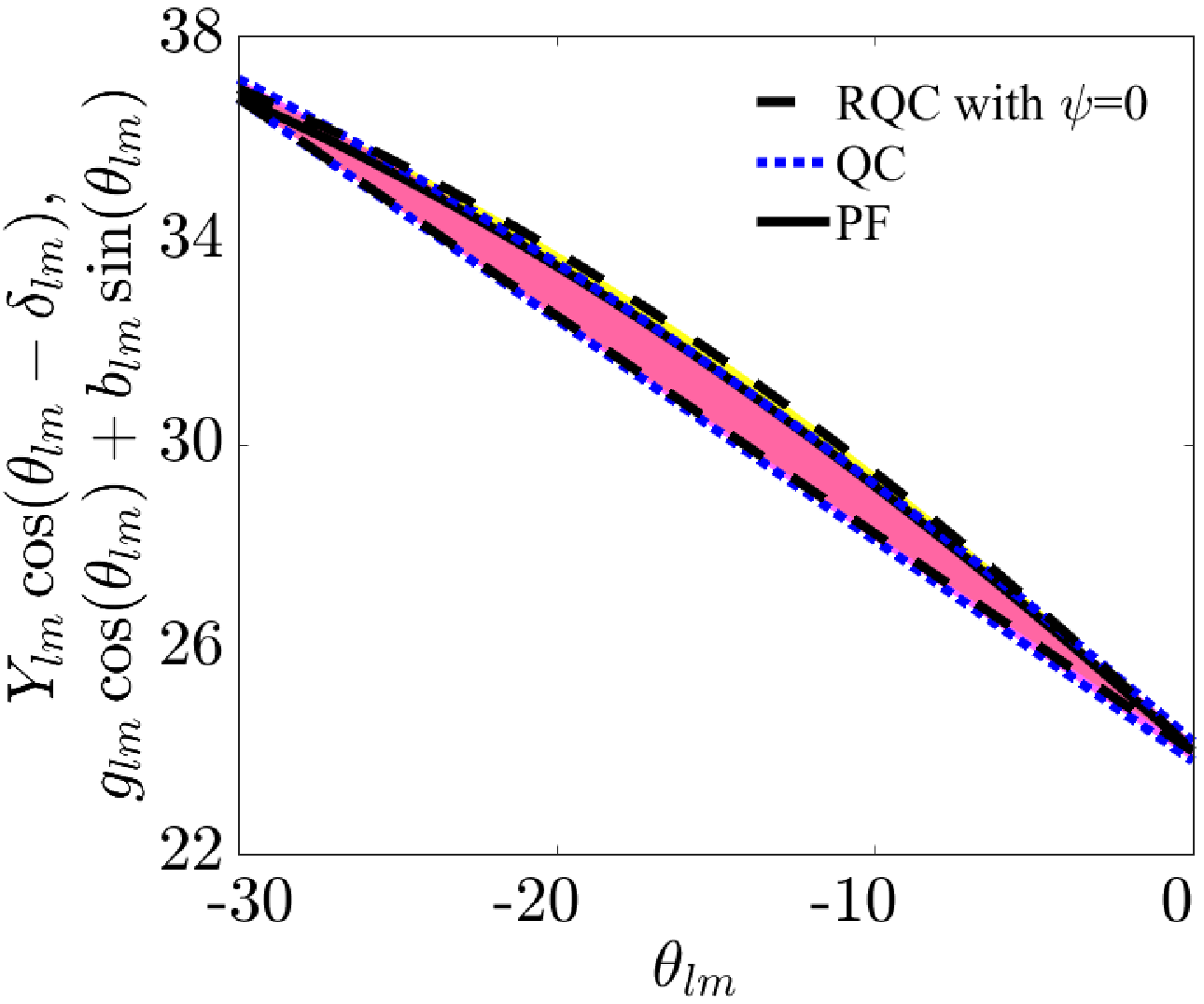}}
\caption{Comparison of envelopes for the trigonometric terms in~\eqref{OPF formulation} and~\eqref{shifted OPF}. The yellow and magenta regions (with dotted and dashed borders, respectively) in (a)--(d) show the envelopes $g_{lm}\left\langle\cos\left(\theta_{lm}\right)\right\rangle^C + b_{lm} \left\langle\sin\left(\theta_{lm}\right)\right\rangle^S$ and $Y_{lm}\left\langle\cos\left(\theta_{lm}-\delta_{lm}\right)\right\rangle^C$, respectively. The black solid lines correspond to the function $g_{lm} \cos\left(\theta_{lm}\right) + b_{lm} \sin\left(\theta_{lm}\right) = Y_{lm} \cos\left(\theta_{lm}-\delta_{lm}\right)$.}
\label{envelops}
\vspace{2em}
\end{figure}

Fig.~\ref{envelops} on the following page illustrates examples of these envelopes for a line with the same mutual admittance ($g_{lm} + jb_{lm} = 0.6 - j 0.8$) for different intervals of angle differences ($\underline{\theta}_{lm} \leq \theta_{lm} \leq \overline{\theta}_{lm}$). While transmission lines with such large resistances are atypical, we choose this admittance to better visualize our approach in Fig.~\ref{envelops}. Our approach is valid for all values of line admittances.


To compare these envelopes, we consider their boundaries.
As shown in%
\ifarxiv
~Appendix~\ref{appendix:Proof of Tighter Lower Bound of RQC relaxation},
\fi either the upper or lower boundary of the envelope $Y_{lm}\left\langle\cos\left(\theta_{lm}-\delta_{lm}\right)\right\rangle^C$ is at least as tight (and sometimes tighter) compared to the corresponding boundary of the envelope $g_{lm}\left\langle\cos\left(\theta_{lm}\right)\right\rangle^C + b_{lm}\left\langle\sin\left(\theta_{lm}\right)\right\rangle^S$ for certain values of $\delta_{lm}$, $\theta_{lm}^{min}$, and $\theta_{lm}^{max}$. In this case, there is no general dominance relationship for the other boundary.
%
%
For other values of $\delta_{lm}$, $\theta_{lm}^{min}$, and $\theta_{lm}^{max}$, none of the boundaries of the envelope $Y_{lm}\left\langle\cos\left(\theta_{lm}-\delta_{lm}\right)\right\rangle^C$ dominate or are dominated by a boundary of the envelope $g_{lm}\left\langle\cos\left(\theta_{lm}\right)\right\rangle^C + b_{lm}\left\langle\sin\left(\theta_{lm}\right)\right\rangle^S$. 
Thus, a QC relaxation that enforces the intersection of these envelopes is generally tighter than relaxations constructed using either of these envelopes individually. Section~\ref{Tightned QC relaxation for Rotated base OPF} discusses this further. 

\begin{figure}[t!]
\centering
\subfloat[$\left\langle\cos\left(\theta_{lm}-\delta_{lm}-\psi\right)\right\rangle^C$]{\includegraphics[width=7.5cm]{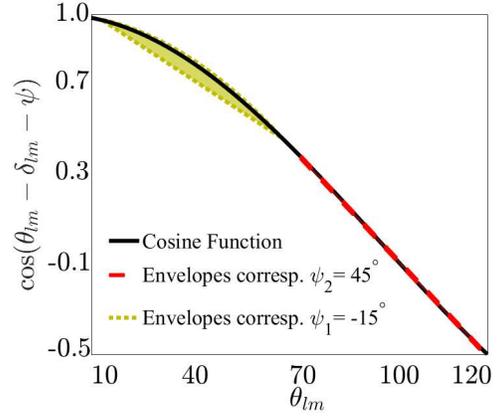}} \\
\subfloat[$\left\langle\sin\left(\theta_{lm}-\delta_{lm}-\psi\right)\right\rangle^S$]{\includegraphics[width=7.5cm]{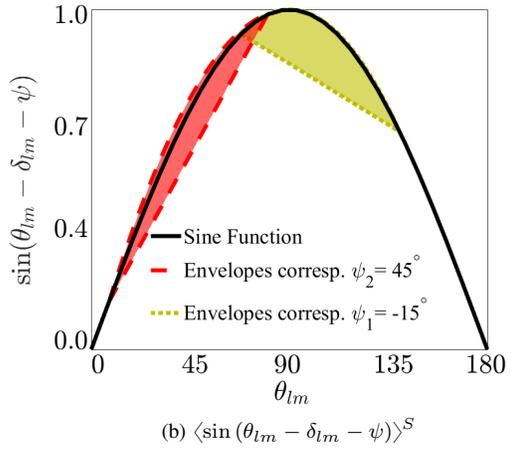}}
\centering
\caption{Comparison of envelopes for the sine and cosine functions for different values of $\psi$. The yellow and red regions (with dashed and dotted borders, respectively) in (a) and (b) show the envelopes $\left\langle\cos\left(\theta_{lm}-\delta_{lm}-\psi\right)\right\rangle^C$ and $\left\langle\sin\left(\theta_{lm}-\delta_{lm}-\psi\right)\right\rangle^S$, for $\psi_1=-15^\circ$ and $\psi_2=45^\circ$, respectively. The angle difference $\theta_{lm}$ varies within $0^\circ\leq\theta_{lm} \leq 72^\circ$, and $\delta_{lm}=-53^\circ$.}
\label{fig:psi_impact}
\end{figure}

The second characteristic distinguishing between the envelopes for~\eqref{OPF formulation} and~\eqref{shifted OPF} is the ability to choose $\psi$ in the latter envelopes. As shown in Fig.~\ref{fig:psi_impact}, changing $\psi$ rotates the arguments of these envelopes.
\mrn{We also visualize the impacts that different values of $\psi$ have on the sine and cosine envelopes in an animation available  at~\url{https://arxiv.org/src/1912.05061v3/anc/rotated_envelope_animation.gif}.}

Analytically comparing the impacts of different values for $\psi$ is not straightforward. Accordingly, this section will later describe an empirical study that suggests a good choice for $\psi$ for typical OPF problems.

\subsection{Envelopes for Trilinear Terms}
\label{Envelopes for trilinear terms}


\mrn{In addition to the trigonometric functions considered thus far, the products between the voltage magnitudes and the trigonometric functions in~\eqref{rotated power flow} are another source of non-convexity in the rotated OPF problem~\eqref{shifted OPF}. We next exploit the relationship between the sending and receiving end trigonometric functions~\eqref{eq:terminal_relationships} in order to relax these products using a limited number of additional lifted variables.}


Similar to~\eqref{eq:qc_cik}--\eqref{eq:qc_sik}, we relax the trilinear products by constructing linear envelopes using the upper and lower bounds on $V_l$, $V_m$, $\cos(\theta_{lm}-\delta_{lm}-\psi)$, $\sin(\theta_{lm}-\delta_{lm}-\psi)$, $\cos(\theta_{lm}+\delta_{lm}+\psi)$, and $\sin(\theta_{lm}+\delta_{lm}+\psi)$. We use the linear relationship~\eqref{eq:terminal_relationships} to represent the upper and lower bounds on the receiving end quantities $\cos(\theta_{lm}+\delta_{lm}+\psi)$ (denoted $\vphantom{\tilde{C}_{lm}}\smash{\underline{\tilde{C}}}_{lm}^{(r)}$, $\vphantom{\tilde{C}_{lm}}\smash{\overline{\tilde{C}}}_{lm}^{(r)}$) and $\sin(\theta_{lm}+\delta_{lm}+\psi)$ (denoted $\vphantom{\tilde{S_{lm}}}\smash{\underline{\tilde{S}}}_{lm}^{(r)}$, $\vphantom{\tilde{S_{lm}}}\smash{\overline{\tilde{S}}}_{lm}^{(r)}$) in terms of the bounds on the sending end quantities $\cos(\theta_{lm}-\delta_{lm}-\psi)$ (denoted $\vphantom{\tilde{C}_{lm}}\smash{\underline{\tilde{C}}}_{lm}^{(s)}$, $\vphantom{\tilde{C}_{lm}}\smash{\overline{\tilde{C}}}_{lm}^{(s)}$) and $\sin(\theta_{lm}-\delta_{lm}-\psi)$ (denoted $\vphantom{\tilde{S}_{lm}}\smash{\underline{\tilde{S}}}_{lm}^{(s)}$, $\vphantom{\tilde{S}_{lm}}\smash{\overline{\tilde{S}}}_{lm}^{(s)}$). We then enforce constraints on the sending end quantities derived from the intersection of the transformed bounds associated with the receiving end quantities along with the bounds on the sending end quantities. Intersecting these bounds forms a polytope in terms of the sending end quantities $\tilde{C}_{lm}^{(s)}$ (representing $\cos(\theta_{lm}-\delta_{lm}-\psi)$) and $\tilde{S}_{lm}^{(s)}$ (representing $\sin(\theta_{lm}-\delta_{lm}-\psi)$), expressible as a convex combination of its extreme points. 

\begin{figure}[t!]
    \centering
\begin{tikzpicture}
\node at (0,0) {\includegraphics[scale=0.5]{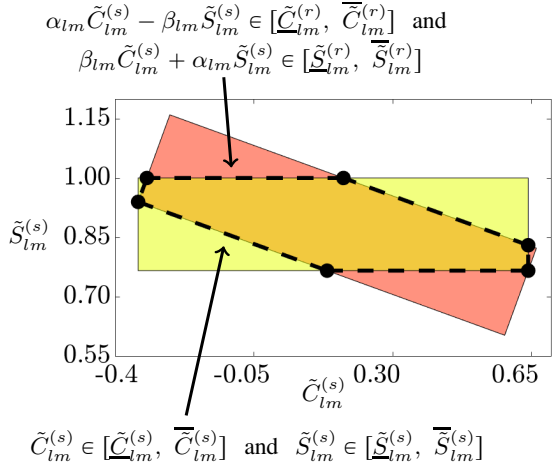}};
\node [anchor=west, font=\small] at (-3.7,2.5) {\begin{tabular}{c} $\alpha_{lm}\tilde{C}_{lm}^{(s)}-\beta_{lm}\tilde{S}_{lm}^{(s)} \in [\vphantom{\tilde{C}_{lm}}\smash{\underline{\tilde{C}}}_{lm}^{(r)},~\vphantom{\tilde{C}_{lm}}\smash{\overline{\tilde{C}}}_{lm}^{(r)}]$ $\text{~and~}$\\[0.3em] $\beta_{lm}\tilde{C}_{lm}^{(s)}+\alpha_{lm}\tilde{S}_{lm}^{(s)} \in [\vphantom{\tilde{S}_{lm}}\smash{\underline{\tilde{S}}}_{lm}^{(r)},~\vphantom{\tilde{S}_{lm}}\smash{\overline{\tilde{S}}}_{lm}^{(r)}]$ \end{tabular}};
\node [anchor=west, font=\small] at (-3.7,-2.9) {$\tilde{C}_{lm}^{(s)} \in [\vphantom{\tilde{C}_{lm}}\smash{\underline{\tilde{C}}}_{lm}^{(s)},~\vphantom{\tilde{C}}\smash{\overline{\tilde{C}}}_{lm}^{(s)}]$ $ \text{~and~}$ $\tilde{S}_{lm}^{(s)} \in [\vphantom{\tilde{S}_{lm}}\smash{\underline{\tilde{S}}}_{lm}^{(s)},~\vphantom{\tilde{S}_{lm}}\smash{\overline{\tilde{S}}}_{lm}^{(s)}]$};
\node [anchor=west, font=\small] at (-0.3,-2.2) {\begin{tabular}{c} $\tilde{C}_{lm}^{(s)}$ \end{tabular}};
\node [anchor=west, font=\small] at (-4.2,-0.1) {\begin{tabular}{c} $\tilde{S}_{lm}^{(s)}$ \end{tabular}};
     	        \draw [-> , line width=0.4mm, black] (-1.5,-2.5) -- (-1.0,-.25);
        \draw [ -> , line width=0.4mm, black] (-1.0,2.0) -- (-0.9,0.8);
\end{tikzpicture}%
	\caption{A projection of the four-dimensional polytope associated with the trilinear products between voltage magnitudes and trigonometric functions, in terms of the sending end variables $\tilde{S}_{lm}^{(s)}$ and $\tilde{C}_{lm}^{(s)}$ representing $\cos(\theta_{lm}-\delta_{lm}-\psi)$ and $\sin(\theta_{lm}-\delta_{lm}-\psi)$. The polytope formed by intersecting the sending end polytope (yellow) and receiving end polytope (red) is outlined with the dashed black lines and has vertices shown by the black dots.}
	\label{fig_intersection}
\end{figure}

Fig.~\ref{fig_intersection} on the following page shows the bounds on both the sending and receiving end quantities in terms of the sending end quantities. The yellow region is the polytope formed by the bounds on $\cos(\theta_{lm}-\delta_{lm}-\psi)$ and $\sin(\theta_{lm}-\delta_{lm}-\psi)$. The red region is the polytope formed by using~\eqref{eq:terminal_relationships} to represent the bounds on the receiving end quantities $\cos(\theta_{lm}+\delta_{lm}+\psi)$ and $\sin(\theta_{lm}+\delta_{lm}+\psi)$ in terms of the sending end quantities $\cos(\theta_{lm}-\delta_{lm}-\psi)$ and $\sin(\theta_{lm}-\delta_{lm}-\psi)$. The black dots are the vertices of the polytope shown by the dashed black lines formed from the intersection of the yellow and red polytopes. 
\ifarxiv
Appendix~\ref{apx:intersection} shows how to compute these vertices.
\fi

Enforcing the constraints associated with both the yellow and red polytopes adds an unnecessary computational burden. We instead restrict the sending end quantities $\cos(\theta_{lm}-\delta_{lm}-\psi)$ and $\sin(\theta_{lm}-\delta_{lm}-\psi)$ to lie within the polytope shown by the black dashed line in Fig.~\ref{fig_intersection}. This implicitly ensures satisfaction of the bounds on the receiving end quantities. 

To relax the product terms $V_l V_m \cos(\theta_{lm}-\delta_{lm}-\psi)$ and $V_l V_m \sin(\theta_{lm}-\delta_{lm}-\psi)$, we first represent the quantities $\cos(\theta_{lm}-\delta_{lm}-\psi)$ and $\sin(\theta_{lm}-\delta_{lm}-\psi)$ using lifted variables $\tilde{C}_{lm}^{(s)}$ and $\tilde{S}_{lm}^{(s)}$, respectively. We then extend the polytope shown by the black dashed lines in Fig.~\ref{fig_intersection} using the upper and lower bounds on $V_l$ and $V_m$. The resulting four-dimensional polytope is the convex hull of the quadrilinear polynomial $V_lV_m\tilde{C}_{lm}^{(s)}\tilde{S}_{lm}^{(s)}$, which we represent using an extreme point formulation similar to~\eqref{eq:qc_cik}--\eqref{eq:qc_sik}. Let $\mathcal{T}_{lm} = \{(\tilde{C}_{lm}^{int,1},\tilde{S}_{lm}^{int,1}),\, (\tilde{C}_{lm}^{int,2},\tilde{S}_{lm}^{int,2}),\,\ldots,\, (\tilde{C}_{lm}^{int,\tilde{N}},\tilde{S}_{lm}^{int,\tilde{N}})\}$ denote the coordinates of the intersection points (black dots) in Fig.~\ref{fig_intersection}, where $\tilde{N}$ is the number of intersection points which ranges from $4$ to $8$ depending on the value of $\psi$. The extreme points of $V_lV_m\tilde{C}_{lm}^{(s)}\tilde{S}_{lm}^{(s)}$ are then denoted as $\eta^{(k)} \in [\underline{V_l},\overline{V_l}] \times [\underline{V_m},\overline{V_m}]\times\mathcal{T}_{lm}$, $k=1,\ldots,4\tilde{N}$. The auxiliary variables $\lambda_k\in[0,1]$, $k=1,\ldots,4\tilde{N}$, are used to form the convex hull of the quadrilinear term $V_lV_m\tilde{C}_{lm}^{(s)}\tilde{S}_{lm}^{(s)}$.

The envelopes for the trilinear terms are:
\begin{align}
    \nonumber & \tilde{c}_{lm} = \!\!\!\! \sum_{k=1,\ldots,4\tilde{N}}\!\! \lambda_k\, \eta^{(k)}_1 \eta^{(k)}_2 \eta^{(k)}_3\!\!, \;\;\;\; \tilde{s}_{lm} = \!\!\!\! \sum_{k=1,\ldots,4\tilde{N}}\!\!\! \lambda_k\, \eta^{(k)}_1 \eta^{(k)}_2 \eta^{(k)}_4\!\!, \\
    \nonumber & V_l =\!\!\!\!\! \sum_{k=1,\ldots,4\tilde{N}} \!\!\!\!\!\lambda_k\eta^{(k)}_1,\;\;
    V_m =\!\!\!\!\! \sum_{k=1,\ldots,4\tilde{N}}\!\!\!\!\! \lambda_k\eta^{(k)}_2\!\!,\;\; 
    \tilde{S}_{lm}^{(s)} =\!\!\!\!\! \sum_{k=1,\ldots,4\tilde{N}}\!\!\!\!\! \lambda_k\eta^{(k)}_4\, ,\\
    \nonumber & \tilde{C}_{lm}^{(s)} =\!\!\!\!\!\!\!\!\!\! \sum_{k=1,\ldots,4\tilde{N}} \!\!\!\!\!\!\!\lambda_k\eta^{(k)}_3\, \!\!\!,
     \!\!\! \sum_{k=1,\ldots,4\tilde{N}} \!\!\!\!\lambda_k = 1, \qquad \!\!\!\!\!\!\lambda_k \geqslant 0, \quad \!\!\!k=1,\ldots,4\tilde{N},\\
      \nonumber & \tilde{C}_{lm}^{(s)} \in \left\langle\cos(\theta_{lm}-\delta_{lm}-\psi)\right\rangle^{C},\\ 
      \nonumber & \tilde{S}_{lm}^{(s)} \in\left\langle\sin(\theta_{lm}-\delta_{lm}-\psi)\right\rangle^{S},\\
     \nonumber &\alpha_{lm} \tilde{S}_{lm}^{(s)} + \beta_{lm} \tilde{C}_{lm}^{(s)} \in \left\langle\sin(\theta_{lm}+\delta_{lm}+\psi)\right\rangle^{S}, \\
    \label{eq:convexhull trilinear terms} &-\beta_{lm} \tilde{S}_{lm}^{(s)} + \alpha_{lm} \tilde{C}_{lm}^{(s)} \in \left\langle\cos(\theta_{lm}+\delta_{lm}+\psi)\right\rangle^{C},
\end{align}
where the last four trigonometric envelope constraints correspond to~\eqref{eq:rotated_envelopes_sending}--\eqref{eq:rotated_envelopes_receiving}.
Note that~\eqref{eq:convexhull trilinear terms} precludes the need for the linking constraint~\mrn{\eqref{eq:qc_linking}} that relates the common term $V_l V_m$ in the products $V_l V_m \sin(\theta_{lm})$ and $V_l V_m \cos(\theta_{lm})$.

\subsection{QC Relaxation of the Rotated OPF Problem}
\label{QC relaxation for Rotated base OPF}
Replacing the squared and trilinear terms with the corresponding lifted variables in the rotated OPF formulation~\eqref{shifted OPF} results in our proposed ``Rotated QC'' (RQC) relaxation: \pagebreak[2]
\begin{subequations}
\label{eq:RQC OPF}
\begin{small}
\begin{align}
& \min\quad \eqref{eq:shifted OPF objective function}\\
&\nonumber \text{subject to} \quad \left(\forall i\in\mathcal{N}, \forall   \left(l,m\right) \in\mathcal{L}\right)\\
& \nonumber \tilde{P}_i^g-\tilde{P}_i^d =\left(g_{sh,i}  \cos(\psi)-b_{sh,i} \sin(\psi)\right)w_{ii}\\ \label{eq:RQC active power injection}&\quad\quad\quad\quad +\sum_{\footnotesize{\substack{(l,m)\in \mathcal{L},\\ \text{s.t.} \hspace{3pt} l=i}}} \tilde{P}_{lm}+\sum_{\footnotesize{\substack{(l,m)\in \mathcal{L},\\ \text{s.t.} \hspace{3pt} m=i}}} \tilde{P}_{ml},\\
& \nonumber \tilde{Q}_i^g-\tilde{Q}_i^d =-\left(g_{sh,i}\sin(\psi)+b_{sh,i}  \cos(\psi)\right)w_{ii}\\ \label{eq:RQC reactive power injection} &\quad\quad\quad\quad +\sum_{\footnotesize{\substack{(l,m)\in \mathcal{L},\\ \text{s.t.} \hspace{3pt} l=i}}} \tilde{Q}_{lm}+\sum_{\footnotesize{\substack{(l,m)\in \mathcal{L},\\ \text{s.t.} \hspace{3pt} m=i}}} \tilde{Q}_{ml},\\
\label{eq:RQC active power flow lm}
& \tilde{P}_{lm} = \left(Y_{lm}\cos(\delta_{lm} + \psi) - b_{c,lm}/2 \sin(\psi)\right)w_{ll}  - Y_{lm}\tilde{c}_{lm},\\[-0.35em]
  \label{eq:RQC reactive power flow lm}
& \tilde{Q}_{lm}= -\left(Y_{lm}\sin(\delta_{lm}+\psi)+b_{c,lm}/2\cos(\psi)\right)w_{ll}-
  Y_{lm}\tilde{s}_{lm},\\[-0.35em]
  \label{eq:RQC active power flow ml}
& \tilde{P}_{ml}=
  -Y_{lm}\tilde{c}_{lm} + \left( Y_{lm} \cos(\delta_{lm}+ \psi) - b_{c,lm}/2 \sin(\psi)\right)w_{mm},\\[-0.35em]
  \label{eq:RQC reactive power flow ml}
& \tilde{Q}_{ml} = 
  Y_{lm}\tilde{s}_{lm} - \left(Y_{lm}  \sin(\delta_{lm}+\psi) + b_{c,lm}/2 \cos(\psi)\right)w_{mm},\\
 \label{eq:RQC squared current}
& \tilde{P}_{lm}^2 + \tilde{Q}_{lm}^2 \leq w_{ll}\, \tilde{\ell}_{lm},\\
 \label{eq:RQC squared current expression}
&\nonumber \tilde{\ell}_{lm} =\bigg(b_{c,lm}^2/4+Y_{lm}^2- Y_{lm}b_{c,lm}\cos(\delta_{lm}+\psi)\sin(\psi)\\[-0.35em]
&\nonumber\quad\quad+Y_{lm}b_{c,lm}\sin(\delta_{lm}+\psi)\cos(\psi)\bigg)\mrn{w_{ll}}+Y_{lm}^2 \mrn{w_{mm}}\\&\quad\quad\nonumber+\left(-2Y_{lm}^2\cos(\delta_{lm}+\psi)+Y_{lm}b_{c,lm}\sin(\psi)\right)\tilde{c}_{lm}\\&\quad\quad
+\left(2Y_{lm}^2\sin(\delta_{lm}+\psi)+Y_{lm}b_{c,lm}\cos(\psi)\right)\tilde{s}_{lm},\\
&  \text{Equations~}\eqref{eq:qc_V},\;\eqref{eq:shifted OPF slack bus}\text{--}\eqref{eq:shifted OPF power flow limit lm},\;\eqref{eq:convexhull trilinear terms}.
\end{align}%
\end{small}%
\end{subequations}%
Note that trilinear terms are relaxed via the extreme point approach in~\eqref{eq:convexhull trilinear terms} that yields the convex hulls for these terms.
The variables $\tilde{c}_{lm}$ and $\tilde{s}_{lm}$ are relaxations of the trilinear terms $V_lV_m\cos(\theta_{lm}-\delta_{lm}-\psi)$ and $V_lV_m\sin(\theta_{lm}-\delta_{lm}-\psi)$, respectively.
\ifarxiv
Appendix~\ref{appendix:power flow formulation with transformer}
\fi
gives an expression for $\tilde{\ell}_{lm}$ that considers off-nominal tap ratios and non-zero phase shifts.

%
%
%
\subsection{Tightened QC Relaxation of the Rotated OPF Problem}
\label{Tightned QC relaxation for Rotated base OPF}

Applying the angle sum and difference identities in combination with~\eqref{eq:terminal_relationships} reveals a linear relationship between the trigonometric functions used in the original QC relaxation~\eqref{eq:QC relaxation}, $\cos(\theta_{lm})$ and $\sin(\theta_{lm})$, and those in the RQC relaxation~\eqref{eq:RQC OPF}, $\cos(\theta_{lm} - \delta_{lm} - \psi)$ and $\sin(\theta_{lm} - \delta_{lm} - \psi)$:
\begingroup \small
\begin{align}
\label{eq:angle_sumdiff3}
& \begin{bmatrix}
\cos(\theta_{lm})\\
\sin(\theta_{lm}) 
\end{bmatrix} = M_{lm}
\begin{bmatrix} 
\sin(\theta_{lm} - \delta_{lm} - \psi)\\ 
\cos(\theta_{lm} - \delta_{lm} - \psi)
\end{bmatrix},
\end{align}\endgroup
\noindent where the constant matrix $M_{lm}$ is defined as
\begin{align*}
M_{lm} = & \frac{1}{2}\left(
\begin{bmatrix}
\sin(\delta_{lm} + \psi) & \cos(\delta_{lm} + \psi) \\
\cos(\delta_{lm} + \psi) & -\sin(\delta_{lm} + \psi)
\end{bmatrix}\begin{bmatrix}
\alpha_{lm} & \beta_{lm} \\
-\beta_{lm} & \alpha_{lm}
\end{bmatrix} \right. \\ 
 & \qquad \left. + \begin{bmatrix}
-\sin(\delta_{lm} + \psi) & \cos(\delta_{lm} + \psi) \\
\cos(\delta_{lm} + \psi) & \sin(\delta_{lm} + \psi)
\end{bmatrix}\right)%
\end{align*}%
with $\alpha_{lm}$ and $\beta_{lm}$ defined as in~\eqref{eq:terminal_relationships}. As mentioned in Section~\ref{Different envelope epresentation explanations}, the RQC relaxation~\eqref{eq:RQC OPF} can be further tightened by additionally enforcing the envelopes $\left\langle\cos\left(\theta_{lm}\right)\right\rangle^C$ and $\left\langle\sin\left(\theta_{lm}\right)\right\rangle^S$ used in the original QC relaxation~\eqref{eq:QC relaxation}. 
This results in the ``Tightened Rotated QC'' (TRQC) relaxation:
%
\begin{subequations}
\label{eq:TRQC OPF}
\small
\begin{align}
& \min\quad \eqref{eq:shifted OPF objective function}\\
&\nonumber \text{subject to} \quad \left(\forall i\in\mathcal{N}, \forall   \left(l,m\right) \in\mathcal{L}\right)\\
& \quad M_{lm}\begin{bmatrix}
\tilde{C}^{(s)}_{lm} \\
\tilde{S}^{(s)}_{lm}
\end{bmatrix}\in \begin{bmatrix}
\left\langle \cos(\theta_{lm}) \right\rangle^C \\
\left\langle \sin(\theta_{lm}) \right\rangle^S 
\end{bmatrix}\\
&\quad \text{Equations~}\eqref{eq:qc_V},\;\eqref{eq:shifted OPF slack bus}\text{--}\eqref{eq:shifted OPF power flow limit lm},\;\eqref{eq:convexhull trilinear terms},\;\eqref{eq:RQC active power injection}\text{--}\eqref{eq:RQC squared current expression}. 
\end{align}
\end{subequations}

\subsection{An Empirical Analysis for Determining the Rotation $\psi$}
\label{Determining best rotation angle}
The key parameter in our proposed QC formulation is the rotation~$\psi$. \mrn{Choosing an appropriate value for $\psi$ using an analytical method is challenging because, as shown in Fig.~\ref{envelops}, $\psi$ simultaneously affects the envelopes for both the sine and cosine functions such that values of $\psi$ which lead to tighter envelopes for the cosine function can result in looser envelopes for the sine function (and vice-versa). Moreover, the single value for $\psi$ applied to the entire system requires balancing the impacts of $\psi$ among all lines simultaneously. Thus, choosing an appropriate value for $\psi$ is not straightforward. We therefore use the following empirical analysis to choose a value for $\psi$ that works well for a range of test cases.} In our results, we denote the best value of $\psi$ for each case as $\psi^*$.

\begin{figure}[h!]
\centering
    \includegraphics[scale=0.50]{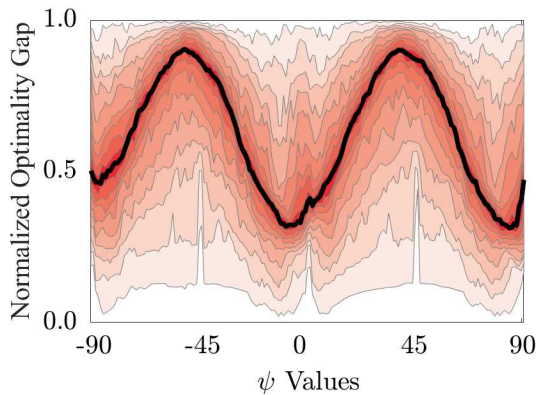}
    \vspace{-.3cm}
    \caption{Normalized optimality gap as a function of $\psi$ for PGLib-OPF cases.}
    \label{normalized_optimality_gap}
\end{figure}

Fig.~\ref{normalized_optimality_gap} shows the optimality gaps for the PGLib-OPF test cases as a function of $\psi$, each normalized by the maximum gap for that case over all values for $\psi$. The results in the figure were generated by sweeping $\psi$ from $-90^\circ$ to $90^\circ$ in steps of $0.5^\circ$. (The figure is exactly symmetric for values of $\psi$ from $90^\circ$ to -$90^\circ$.) The shaded red bands around the median line (in black) show every fifth percentile of the results. 

The results in Fig.~\ref{normalized_optimality_gap} indicate that good values of $\psi$ are consistent across the test systems. Thus, we suggest using $\psi = 80^\circ$, which is where the median of the optimality gaps over all the test cases was smallest. Moreover, the symmetry in Fig.~\ref{normalized_optimality_gap} implies that
selecting $\psi$ within the intervals [$-90^\circ$, $-80^\circ$], [$-15^\circ$, $-5^\circ$], and [$80^\circ$, $90^\circ$] results in nearly the smallest optimality gaps for almost all of the test cases compared to the optimality gaps from the RQC relaxation using $\psi^*$.



\mrn{Our ongoing work is generalizing the rotated QC formulation proposed in this paper by allowing for different values of $\psi$ for each line in the system in order to construct tighter envelopes (i.e., independent choices of $\psi_{lm}$ for each $(l,m)\in\mathcal{L}$ rather than one value of $\psi$ for the entire system). While complicating the overall formulation, this generalization simplifies the impacts from choosing each $\psi_{lm}$ and thus has the potential to enable analytical methods for identifying the best values for $\psi_{lm}$ for each $(l,m)\in\mathcal{L}$. One possible approach is to compute the value of $\psi_{lm}$ that minimizes the areas of the envelopes (i.e., the yellow and red regions in Fig.~\ref{fig:psi_impact}).}



\section{Numerical results}
\label{Numerical_results}
This section demonstrates the effectiveness of the proposed approach using selected test cases from the PGLib-OPF v18.08 benchmark library~\cite{pglib}. 
These test cases were selected since existing relaxations fail to provide tight bounds on the best known objective values.
Our implementations use Julia 0.6.4, JuMP v0.18~\cite{JuMP}, PowerModels.jl~\cite{powermodels}, and Gurobi 8.0 as modeling tools and the solver. The results are computed using a laptop with an i7 1.80~GHz processor and 16 GB of RAM. 

Table~\ref{tab:various_QC_results} summarizes the results from applying the QC~\eqref{eq:QC relaxation}, RQC ~\eqref{eq:RQC OPF}, and TRQC~\eqref{eq:TRQC OPF} relaxations to selected test cases. The first column lists the test cases. The next group of columns represents optimality gaps, defined as
\begin{align}
\label{eq:optimality gap}
\text{\emph{Optimality~Gap}}=\left(\dfrac{\text{\emph{Local~Solution}} - \text{\emph{QC Bound}}}{\text{\emph{Local~Solution}}}\right).
\end{align} 
The optimality gaps are defined using the local solutions to the non-convex problem~\eqref{OPF formulation} from \mbox{PowerModels.jl}. The final group of columns shows the solver times. 

\begin{table*}[t!]
\caption{Results from Applying the QC and RQC Relaxations with Various Options to Selected PGLib Test Cases}
\label{tab:various_QC_results}
\vspace{-0.7cm}
\begin{center}
\begin{tabular}{|p{2.45cm}|c|c|c|c|c|c|c|c||c|c|c|}
\hline 
\multirow{2}{*}{\thead{Test Cases}} & 
\multirow{2}{*}{\thead{QC \\Gap (\%)}} & 
\multirow{2}{*}{\thead{RQC \\ ($\psi=0$) \\ Gap (\%)}} & 
\multirow{2}{*}{\thead{RQC\\($\psi=80^\circ$)\\ Gap (\%)}} & 
\multicolumn{2}{c|}{{\thead{RQC $(\psi^*)$ }}} & 
\multirow{2}{*}{\thead{TRQC \\ ($\psi=80^\circ$)\\ Gap (\%)}} & 
\multicolumn{2}{c|}{{\thead{TRQC ($\psi^*$)}}} & 
\multirow{2}{*}{\thead{QC\\Time\\(sec)}} & 
\multirow{2}{*}{\thead{RQC\\Time\\(sec)}} & 
\multirow{2}{*}{\thead{TRQC\\Time\\(sec)}}\tabularnewline
\cline{5-6} \cline{6-6} \cline{8-9} \cline{9-9} 
   &  & & &\thead{Gap (\%)} & \thead{$\psi^\ast$}& &\thead{Gap (\%)} & \thead{$\psi^\ast$}  &&& \tabularnewline
\hline
case3\_lmbd &  \hphantom{0}0.97 & \hphantom{0}0.97 & \hphantom{0}0.89& \hphantom{0}0.79 & $-81^\circ$ & \hphantom{0}0.84&\hphantom{0}0.63 & \hphantom{00}$11^\circ$ 
& \hphantom{00}0.34&\hphantom{00}0.01 &\hphantom{00}0.01\tabularnewline
\hline 
case30\_ieee &  18.67 & 14.91& 13.14 & 12.11 & \hphantom{00}$65^\circ$&13.14 & 11.82 & $-25^\circ$ 
&\hphantom{00}0.33&\hphantom{00}0.03& \hphantom{00}0.03\tabularnewline
\hline
case118\_ieee &  \hphantom{0}0.77 & \hphantom{0}0.90&\hphantom{0}0.65 & \hphantom{0}0.64 & \hphantom{00}$70^\circ$& \hphantom{0}0.64 & \hphantom{0}0.62 & \hphantom{00}$70^\circ$ 
&\hphantom{00}0.55&\hphantom{00}0.19&\hphantom{00}0.23\tabularnewline
\hline 
case300\_ieee  & \hphantom{0}2.56 & \hphantom{0}2.58&\hphantom{0}2.43 & \hphantom{0}2.26 & $-13^\circ$& \hphantom{0}2.32 & \hphantom{0}2.24 & $-13^\circ$
&\hphantom{00}1.54&\hphantom{00}1.50&\hphantom{00}3.15\tabularnewline
\hline 
case9241\_pegase &  \hphantom{0}1.71 & \hphantom{0}1.70&\hphantom{0}1.70 & \hphantom{0}1.69 & $-10^\circ$& \hphantom{0}1.70 & \hphantom{0}1.69 & $-10^\circ$ 
&265.39&190.80&297.56\tabularnewline
\hline 
\hline 
case3\_lmbd\_\_api &  \hphantom{0}4.57 & \hphantom{0}4.31&\hphantom{0}4.42 & \hphantom{0}4.28 & \hphantom{000}$2^\circ$& \hphantom{0}4.17 & \hphantom{0}3.93 & $-71^\circ$ 
&\hphantom{00}0.51&\hphantom{00}0.01&\hphantom{00}0.01\tabularnewline
\hline 
case24\_ieee\_rts\_\_api &  11.02 & \hphantom{0}7.83 &7.51 & \hphantom{0}7.24 & \hphantom{00}$79^\circ$&\hphantom{0}7.31& \hphantom{0}6.98 & $-11^\circ$ 
& \hphantom{0}0.71&\hphantom{00}0.03&\hphantom{00}0.04\tabularnewline
\hline 
case39\_epri\_\_api &  \hphantom{0}1.71 & \hphantom{0}1.38&\hphantom{0}1.33 & \hphantom{0}1.33 & $-11^\circ$& \hphantom{0}1.32 & \hphantom{0}1.32 & \hphantom{00}$79^\circ$ 
&\hphantom{00}0.39&\hphantom{00}0.05&\hphantom{00}0.05\tabularnewline
\hline 
case73\_ieee\_rts\_\_api &  \hphantom{0}9.54 & \hphantom{0}8.12&\hphantom{0}7.36 & \hphantom{0}7.36 & $-10^\circ$& \hphantom{0}7.24 & \hphantom{0}7.24 & $-10^\circ$ 
&\hphantom{00}1.00&\hphantom{00}0.29&\hphantom{00}0.37\tabularnewline
\hline 
case118\_ieee\_\_api &  28.67 & 28.03&26.82 & 26.52 & \hphantom{0}$-8^\circ$& 27.11 & 26.38 & \hphantom{0}$-8^\circ$ &\hphantom{00}0.53&\hphantom{00}0.20&\hphantom{00}0.97\tabularnewline
\hline 
case179\_goc\_\_api  & \hphantom{0}5.86 & \hphantom{0}6.01&\hphantom{0}5.57 & \hphantom{0}4.90 & $-81^\circ$& \hphantom{0}4.90 & \hphantom{0}4.06 & $-78^\circ$ 
&\hphantom{00}0.82&\hphantom{00}0.61&\hphantom{00}0.64\tabularnewline
\hline 
\hline 
 case14\_ieee\_\_sad &  19.16 & 21.45&17.89 & 16.30 & \hphantom{00}$77^\circ$& 15.82 & 15.39 & $-12^\circ$ 
&\hphantom{00}0.35&\hphantom{00}0.03&\hphantom{00}0.03\tabularnewline
\hline 
 case24\_ieee\_rts\_\_sad  & \hphantom{0}2.74 & \hphantom{0}2.55&\hphantom{0}2.31 & \hphantom{0}2.19 & \hphantom{00}$78^\circ$& \hphantom{0}2.26 & \hphantom{0}2.12 & $-12^\circ$ 
&\hphantom{00}0.40&\hphantom{00}0.05&\hphantom{00}0.06\tabularnewline
\hline
 case30\_ieee\_\_sad &  \hphantom{0}5.66 & \hphantom{0}5.95&\hphantom{0}4.81 & \hphantom{0}4.59 & $-13^\circ$& \hphantom{0}4.56 & \hphantom{0}4.45 & \hphantom{00}$66^\circ$ 
&\hphantom{00}0.32&\hphantom{00}0.05&\hphantom{00}0.06\tabularnewline
\hline 
 case73\_ieee\_rts\_\_sad &  \hphantom{0}2.37 & \hphantom{0}2.24&\hphantom{0}1.98 & \hphantom{0}1.90 & \hphantom{00}$79^\circ$& \hphantom{0}1.84  & \hphantom{0}1.82 & \hphantom{00}$78^\circ$
&\hphantom{00}0.41&\hphantom{00}0.28&\hphantom{00}0.41\tabularnewline
\hline 
case118\_ieee\_\_sad &  \hphantom{0}6.67 & \hphantom{0}8.10&\hphantom{0}5.45 & \hphantom{0}5.39 & \hphantom{00}$81^\circ$& \hphantom{0}5.45 & \hphantom{0}5.07 & \hphantom{00}$69^\circ$ 
&\hphantom{00}0.58&\hphantom{00}0.25&\hphantom{00}0.39\tabularnewline
\hline 
case162\_ieee\_dtc\_\_sad &  \hphantom{0}6.22 & \hphantom{0}6.30&\hphantom{0}5.65 & \hphantom{0}5.59 & $-14^\circ$& \hphantom{0}5.65 & \hphantom{0}5.54 & \hphantom{00}$76^\circ$ 
&\hphantom{00}0.86&\hphantom{00}0.55&\hphantom{00}0.84\tabularnewline
\hline 
case300\_ieee\_\_sad &  \hphantom{0}2.34 & \hphantom{0}2.59&\hphantom{0}1.80 & \hphantom{0}1.61 & \hphantom{00}$83^\circ$& \hphantom{0}1.78 & \hphantom{0}1.59 & \hphantom{00}$83^\circ$ 
&\hphantom{00}1.94&\hphantom{00}1.29&\hphantom{00}2.06\tabularnewline
\hline 
\end{tabular}
{\raggedright \emph{AC}: AC local solution from~\eqref{OPF formulation}, \emph{QC Gap}: Optimality gap for the QC relaxation from~\eqref{eq:QC relaxation}, \emph{RQC Gap}: Optimality gap for the Rotated QC relaxation from~\eqref{eq:RQC OPF},\newline \emph{TRQC Gap}: Optimality gap for the Tightened Rotated QC Relaxation from~\eqref{eq:TRQC OPF}, $\psi^*$: Use of the best $\psi$ for this case. \par}
\vspace{-0.5cm}
\end{center}
\end{table*}


Comparing the second and third columns in Table~\ref{tab:various_QC_results} reveals that using admittances in polar form without rotation (i.e., the RQC relaxation~\eqref{eq:RQC OPF} with $\psi = 0$) can improve the optimality gaps of some test cases (e.g., improvements of 3.76$\%$ and 3.19$\%$ for ``case30$\_$ieee'' and ``case24$\_$ieee$\_$rts$\_\_$api'', respectively, relative to the original QC relaxation~\eqref{eq:QC relaxation}) .
However, the RQC relaxation with $\psi=0$ has worse performance in other cases, such as ``case300$\_$ieee'' and ``case14$\_$ieee$\_\_$sad'', which have 0.02$\%$ and 2.29$\%$ larger optimality gaps, respectively. 

Using a non-zero value for $\psi$ can improve the optimality gaps. Solving the RQC relaxation~\eqref{eq:RQC OPF} with the suggested $\psi = 80^\circ$ obtained from the empirical analysis in Section~\ref{Determining best rotation angle} results in 1.08\% better optimality gaps, on average, compared to the original QC relaxation. 
The RQC relaxation~\eqref{eq:RQC OPF} with $\psi^*$ (the best value of $\psi$ for each case) provides optimality gaps that are not worse than those obtained by the original QC relaxation~\eqref{eq:QC relaxation} for all test cases, yielding an improvement of 1.36\% on average compared to the original QC relaxation. 
As one specific example, the gap from the original QC relaxation for ``case162$\_$ieee$\_$dtc$\_\_$sad'' is $6.22\%$ compared to $6.30\%$ for the RQC relaxation~\eqref{eq:RQC OPF} with $\psi=0$ relaxation~\eqref{eq:RQC OPF}. Use of the suggested $\psi = 80^\circ$ reduces the gap to 5.65$\%$, which is superior to the gap obtained from the QC relaxation~\eqref{eq:QC relaxation}. Using $\psi^*$ further reduces the optimality gap to 5.59$\%$.

Enforcing the envelopes from both the original QC relaxation and the RQC relaxation, i.e., the TRQC relaxation~\eqref{eq:TRQC OPF}, further improves the optimality gaps. Solving the TRQC relaxation~\eqref{eq:TRQC OPF} with the suggested $\psi = 80^\circ$ results in 1.29\% better gaps, on average, compared to the original QC relaxation. The TRQC relaxation with $\psi^*$ yields optimality gaps that are $1.57$\% and $0.21$\% better, on average, compared to the original QC relaxation and the RQC relaxation with $\psi^*$. The additional envelopes $\left\langle\sin\left(\theta_{lm} \right)\right\rangle^S$ and $\left\langle\cos\left(\theta_{lm} \right)\right\rangle^C$ in the TRQC relaxation increase the average solver time by $22$\%.

\begin{figure}[t!]
\vspace{-0.3cm}
\centering
    \includegraphics[scale=0.5]{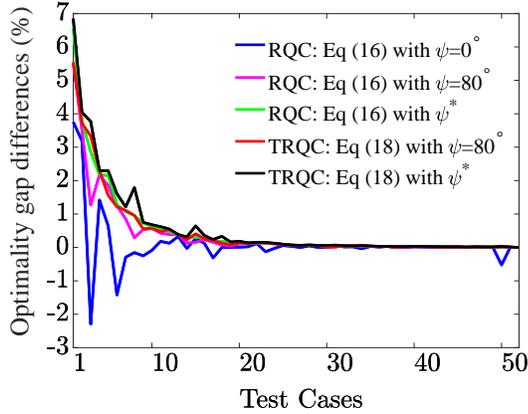}
    \caption{Comparison of optimality gap differences with respect to the original QC relaxation~\eqref{eq:QC relaxation} for different QC relaxation variants.}
    \label{optimality_gap_differnce}
    \vspace*{-1.5em}
\end{figure}

Fig.~\ref{optimality_gap_differnce} visualizes the optimality gaps for variants of the QC relaxation over a range of test cases. Positive values indicate an improvement in the optimality gap of the associated variant relative to the original QC relaxation~\eqref{eq:QC relaxation}. The test cases are sorted in order of increasing optimality gaps obtained from the original QC relaxation. The TRQC relaxation with $\psi^*$ achieves the smallest optimality gaps. While the RQC relaxation with $\psi=0$ obtains a worse optimality gap for some test cases compared to the original QC relaxation, both the RQC and TRQC relaxations with $\psi^*$ outperform the QC relaxation for all test cases. As expected from the analysis in Section~\ref{Determining best rotation angle}, applying the suggested $\psi = 80^\circ$ results in good performance across a variety of test cases.

\section{Conclusion}
\label{conclusion}
This paper proposes and empirically tests two improvements for strengthening QC relaxations of OPF problems by tightening the envelopes used for the trigonometric terms. The first improvement represents the line admittances in polar form. The second improvement applies a complex base power normalization with angle $\psi$ in order to rotate the arguments of the trigonometric terms. An empirical analysis is used to suggest a good value for $\psi$. Comparison to the state-of-the-art QC relaxation reveals the effectiveness of the proposed improvements. Our ongoing work is extending the RQC relaxation to allow for distinct values of $\psi$ for each line.

\ifarxiv
\appendix
\section{Appendices}
\subsection{\mrn{Sine and Cosine Envelopes for Extended Angle Ranges}}
\label{appendix: sine and cosine envelopes}
\mrn{Tables~\ref{tab:sine_envelopes} and~\ref{tab:cosine_envelopes} on the following page extend the envelopes for the sine and cosine functions, $\left\langle \sin(x)\right\rangle^S$ and $\left\langle \cos(x)\right\rangle^c$, defined in~\eqref{eq:sine envelope} and~\eqref{eq:cosine envelope} in order to consider all angle ranges.}

\begin{table}
\begin{center}
\caption{\mrn{Sine Function Envelopes For Extended Angle Ranges}}
\label{tab:sine_envelopes}
\begin{tabular}{|c|c|}
\hline 
Envelope & Interval\tabularnewline
\hline 
\hline 
\thead{$\widecheck{S}\leq 1-\frac{1-\cos\left({x^m}\right)}{\left(x^m\right)^2}\left(x+\frac{3\pi}{2}\right)^2,$\\
$\widecheck{S}\geq\frac{\sin\left(\underline{x}\right)-\sin\left({\overline{x}}\right)}{{\underline{x}-\overline{x}}}\left(x-{\underline{x}}\right)+\sin\left(\underline{x}\right),$\\$x^m=\max(\left|\underline{x}+\frac{3\pi}{2}\right|,\left|\overline{x}+\frac{3\pi}{2}\right|)$} & \thead{$-2\pi\leq\underline{x}$\\$\overline{x}\leq-\pi$}\tabularnewline
\hline 
\thead{$\widecheck{S}\!\leq-\cos\left(\frac{x^m}{2}\right)\!\!\left(x\!+\pi+\!\frac{x^m}{2}\right)\!\!+\!\sin \left(\frac{x^m}{2}\right),$\\
$\widecheck{S}\!\geq-\cos\left(\frac{x^m}{2}\right)\!\!\left(x\!+\pi-\!\frac{x^m}{2}\right)\!\!-\!\sin\left(\frac{x^m}{2}\right),$\\ $x^m=\max(\left|\underline{x}+\pi\right|,\left|\overline{x}+\pi\right|)$} & \thead{$-\frac{3\pi}{2}\leq\underline{x}\leq-\pi$\\$-\pi\leq\overline{x}\leq -\frac{\pi}{2}$}\tabularnewline
\hline 
\thead{$\widecheck{S}\geq -\left(1-\frac{1-\cos\left({x^m}\right)}{\left(x^m\right)^2}\left(x+\frac{\pi}{2}\right)^2\right),$\\
$\widecheck{S}\leq\frac{\sin\left(\underline{x}\right)-\sin\left({\overline{x}}\right)}{{\underline{x}-\overline{x}}}\left(x-{\underline{x}}\right)+\sin\left(\underline{x}\right),$\\$x^m=\max(\left|\underline{x}+\frac{\pi}{2}\right|,\left|\overline{x}+\frac{\pi}{2}\right|)$} & \thead{$-\pi\leq\underline{x}$\\$\overline{x}\leq0$}\tabularnewline
\hline 
\thead{$\widecheck{S}\!\leq\!\cos\left(\frac{x^m}{2}\right)\!\!\left(x\!-\!\frac{x^m}{2}\right)\!\!+\!\sin \left(\frac{x^m}{2}\right),$\\
$\widecheck{S}\!\geq\!\cos\left(\frac{x^m}{2}\right)\!\!\left(x\!+\!\frac{x^m}{2}\right)\!\!-\!\sin\left(\frac{x^m}{2}\right),$\\ $x^m=\max(\left|\underline{x}\right|,\left|\overline{x}\right|)$} & \thead{$-\frac{\pi}{2}\leq\underline{x}\leq0$\\$0\leq\overline{x}\leq\frac{\pi}{2}$}\tabularnewline
\hline 
\thead{$\widecheck{S}\leq 1-\frac{1-\cos\left({x^m}\right)}{\left(x^m\right)^2}\left(x-\frac{\pi}{2}\right)^2,$\\
$\widecheck{S}\geq\frac{\sin\left(\underline{x}\right)-\sin\left({\overline{x}}\right)}{{\underline{x}-\overline{x}}}\left(x-{\underline{x}}\right)+\sin\left(\underline{x}\right),$\\$x^m=\max(\left|\underline{x}-\frac{\pi}{2}\right|,\left|\overline{x}-\frac{\pi}{2}\right|)$} & \thead{$0\leq\underline{x}$\\$\overline{x}\leq\pi$}\tabularnewline
\hline 
\thead{$\widecheck{S}\!\leq-\cos\left(\frac{x^m}{2}\right)\!\!\left(x\!-\pi+\!\frac{x^m}{2}\right)\!\!+\!\sin \left(\frac{x^m}{2}\right),$\\
$\widecheck{S}\!\geq-\cos\left(\frac{x^m}{2}\right)\!\!\left(x\!-\pi-\!\frac{x^m}{2}\right)\!\!-\!\sin\left(\frac{x^m}{2}\right),$\\ $x^m=\max(\left|\underline{x}-\pi\right|,\left|\overline{x}-\pi\right|)$} & \thead{$\frac{\pi}{2}\leq\underline{x}\leq\pi$\\$\pi\leq\overline{x}\leq\frac{3\pi}{2}$}\tabularnewline
\hline 
\thead{$\widecheck{S}\geq -\left(1-\frac{1-\cos\left({x^m}\right)}{\left(x^m\right)^2}\left(x-\frac{3\pi}{2}\right)^2\right),$\\
$\widecheck{S}\leq\frac{\sin\left(\underline{x}\right)-\sin\left({\overline{x}}\right)}{{\underline{x}-\overline{x}}}\left(x-{\underline{x}}\right)+\sin\left(\underline{x}\right),$\\$x^m=\max(\left|\underline{x}-\frac{3\pi}{2}\right|,\left|\overline{x}-\frac{3\pi}{2}\right|)$} & \thead{$\pi\leq\underline{x}$\\$\overline{x}\leq2\pi$}\tabularnewline
\hline 
\end{tabular}
\end{center}
\end{table}

\begin{table}
\begin{center}
\caption{\mrn{Cosine Function Envelopes For Extended Angle Ranges}}
\label{tab:cosine_envelopes}
\begin{tabular}{|c|c|}
\hline 
Envelope & Interval\tabularnewline
\hline 
\hline 
\thead{$\widecheck{C}\!\leq-\cos\left(\frac{x^m}{2}\right)\!\!\left(x\!+\frac{3\pi}{2}+\!\frac{x^m}{2}\right)\!\!+\!\sin \left(\frac{x^m}{2}\right),$\\
$\widecheck{C}\!\geq-\cos\left(\frac{x^m}{2}\right)\!\!\left(x\!+\frac{3\pi}{2}-\!\frac{x^m}{2}\right)\!\!-\!\sin\left(\frac{x^m}{2}\right),$\\ $x^m=\max(\left|\underline{x}+\frac{3\pi}{2}\right|,\left|\overline{x}+\frac{3\pi}{2}\right|)$} & \thead{$-2\pi\leq\underline{x}\leq-\frac{3\pi}{2}$\\$-\frac{3\pi}{2}\overline{x}\leq-\pi$}\tabularnewline
\hline 
\thead{$\widecheck{C}\geq -\left(1-\frac{1-\cos\left({x^m}\right)}{\left(x^m\right)^2}\left(x+\pi\right)^2\right),$\\
$\widecheck{C}\leq\frac{\cos\left(\underline{x}\right)-\cos\left({\overline{x}}\right)}{{\underline{x}-\overline{x}}}\left(x-{\underline{x}}\right)+\cos\left(\underline{x}\right),$\\$x^m=\max(\left|\underline{x}+\pi\right|,\left|\overline{x}+\pi\right|)$} & \thead{$-\frac{3\pi}{2}\leq\underline{x}$\\$\overline{x}\leq-\frac{\pi}{2}$}\tabularnewline
\hline 
\thead{$\widecheck{C}\!\leq\cos\left(\frac{x^m}{2}\right)\!\!\left(x\!+\frac{\pi}{2}-\!\frac{x^m}{2}\right)\!\!+\!\sin \left(\frac{x^m}{2}\right),$\\
$\widecheck{C}\!\geq\cos\left(\frac{x^m}{2}\right)\!\!\left(x\!+\frac{\pi}{2}+\!\frac{x^m}{2}\right)\!\!-\!\sin\left(\frac{x^m}{2}\right),$\\ $x^m=\max(\left|\underline{x}+\frac{\pi}{2}\right|,\left|\overline{x}+\frac{\pi}{2}\right|)$} & \thead{$-\pi\leq\underline{x}\leq-\frac{\pi}{2}$\\$-\frac{\pi}{2}\leq\overline{x}\leq0$}\tabularnewline
\hline 
\thead{$\widecheck{C}\leq \left(1-\frac{1-\cos\left({x^m}\right)}{\left(x^m\right)^2}\left(x+\pi\right)^2\right),$\\
$\widecheck{C}\geq\frac{\cos\left(\underline{x}\right)-\cos\left({\overline{x}}\right)}{{\underline{x}-\overline{x}}}\left(x-{\underline{x}}\right)+\cos\left(\underline{x}\right),$\\$x^m=\max(\left|\underline{x}\right|,\left|\overline{x}\right|)$} & \thead{$-\frac{\pi}{2}\leq\underline{x}$\\$\overline{x}\leq\frac{\pi}{2}$}\tabularnewline
\hline 
\thead{$\widecheck{C}\!\leq-\cos\left(\frac{x^m}{2}\right)\!\!\left(x\!-\frac{\pi}{2}+\!\frac{x^m}{2}\right)\!\!+\!\sin \left(\frac{x^m}{2}\right),$\\
$\widecheck{C}\!\geq-\cos\left(\frac{x^m}{2}\right)\!\!\left(x\!-\frac{\pi}{2}-\!\frac{x^m}{2}\right)\!\!-\!\sin\left(\frac{x^m}{2}\right),$\\ $x^m=\max(\left|\underline{x}-\frac{\pi}{2}\right|,\left|\overline{x}-\frac{\pi}{2}\right|)$} & \thead{$0\leq\underline{x}\leq\frac{\pi}{2}$\\$\frac{\pi}{2}\leq\overline{x}\leq\pi$}\tabularnewline
\hline 
\thead{$\widecheck{C}\geq -\left(1-\frac{1-\cos\left({x^m}\right)}{\left(x^m\right)^2}\left(x-\pi\right)^2\right),$\\
$\widecheck{C}\leq\frac{\cos\left(\underline{x}\right)-\cos\left({\overline{x}}\right)}{{\underline{x}-\overline{x}}}\left(x-{\underline{x}}\right)+\cos\left(\underline{x}\right),$\\ $x^m=\max(\left|\underline{x}-\pi\right|,\left|\overline{x}-\pi\right|)$} & \thead{$\frac{\pi}{2}\leq\underline{x}$\\$\overline{x}\leq\frac{3\pi}{2}$}\tabularnewline
\hline 
\thead{$\widecheck{C}\!\leq\cos\left(\frac{x^m}{2}\right)\!\!\left(x\!-\frac{3\pi}{2}-\!\frac{x^m}{2}\right)\!\!+\!\sin \left(\frac{x^m}{2}\right),$\\
$\widecheck{C}\!\geq\cos\left(\frac{x^m}{2}\right)\!\!\left(x\!-\frac{3\pi}{2}+\!\frac{x^m}{2}\right)\!\!-\!\sin\left(\frac{x^m}{2}\right),$\\$x^m=\max(\left|\underline{x}-\frac{3\pi}{2}\right|,\left|\overline{x}-\frac{3\pi}{2}\right|)$} & \thead{$\pi\leq\underline{x}\leq\frac{3\pi}{2}$\\$\frac{3\pi}{2}\leq\overline{x}\leq2\pi$}\tabularnewline
\hline 
\end{tabular}
\end{center}
\end{table}

\subsection{More General Line Models}
\label{appendix:power flow formulation with transformer}
This appendix extends the paper's results to a line model that considers transformers with a non-zero phase shift $\theta_{lm}^{shift}$ and/or an off-nominal voltage ratio $\tau_{lm}$. With this model, the complex power flows into both terminals of line $\left(l,m\right)\in\mathcal{L}$ are:
\begin{subequations}
\small
\label{complex power flow with trans}
\begin{align}
\label{complex power flow sending with trans}
 &S_{lm}=V_l e^{j\theta_l}
  \!\left[\left(Y_{lm}e^{j\delta_{lm}}+j\frac{b_{c,lm}}{2}\right)\frac{V_l e^{j\theta_l}}{\tau_{lm}^2} - \frac{Y_{lm} e^{j\delta_{lm}} V_m e^{j\theta_m}}{\tau_{lm} e^{-j\theta_{lm}^{shift}}}\!\right]^*\\
  \label{complex power flow receiving with trans}
 &S_{ml}=V_m e^{j\theta_m}
  \!\left[\left(Y_{lm}e^{j\delta_{lm}}+j\frac{b_{c,lm}}{2}\right) V_m e^{j\theta_m}\! -\! \frac{Y_{lm}e^{j\delta_{lm}}V_l e^{j\theta_l}}{\tau_{lm} e^{j\theta_{lm}^{shift}}}\!\right]^*
\end{align}
\end{subequations}
We follow the procedure in Section~\ref{Rotated Base OPF formulation} by applying a complex base power normalization:
\begin{subequations}
\label{rotated power flow sending with trans}
\begin{align}
\label{rotated power flow sending with trans1}
&\nonumber \tilde{S}_{lm}=\frac{S_{lm}}{e^{j\psi}}= \left(\frac{Y_{lm}}{\tau_{lm}^2}e^{-j(\delta_{lm}+\psi)}+\frac{b_{c,lm}}{2}\frac{e^{-j(\frac{\pi}{2}+\psi)}}{\tau_{lm}^2}\right)V_l^2\\&\qquad\qquad\qquad\quad - \frac{Y_{lm}}{\tau_{lm}}V_l V_m e^{j(\theta_{lm}-\delta_{lm}-\theta_{lm}^{shift}-\psi)},\\
  \label{rotated power flow sending with trans2}
&\nonumber \tilde{S}_{ml} = \frac{S_{ml}}{e^{j\psi}} = \left(Y_{lm}  e^{-j(\delta_{lm}+\psi)}+\frac{b_{c,lm}}{2}  e^{-j(\frac{\pi}{2}+\psi)}\right)V_m^2
  \\&\qquad\qquad\qquad\quad -\frac{Y_{lm}}{\tau_{lm}}V_l V_m e^{j(-\theta_{lm}-\delta_{lm}+\theta_{lm}^{shift}-\psi)}.
\end{align}
\end{subequations}
Taking the real and imaginary parts of~\eqref{rotated power flow sending with trans} yields:
\begin{subequations}
\label{eq:general rotated power flow}
\begin{align}
\nonumber  &\!\!\!\tilde{P}_{lm}\! =\! \Re(\tilde{S}_{lm})\! =\! \left(\frac{Y_{lm}}{\tau_{lm}^2}\cos(\delta_{lm}+\psi)-\frac{b_{c,lm}}{2\tau_{lm}^2}\sin(\psi)\right)V_l^2\\ \label{eq:general shifted OPF active power flow lm} &\qquad\qquad -
  \frac{Y_{lm}}{\tau_{lm}}V_l V_m\cos(\theta_{lm}-\delta_{lm}-\theta_{lm}^{shift}-\psi),\\
  \label{eq:general shifted OPF reactive power flow lm}
\nonumber &\!\!\! \tilde{Q}_{lm} \!=\! \Im(\tilde{S}_{lm}) \!=\! \left(-\frac{Y_{lm}}{\tau_{lm}^2}\sin(\delta_{lm}+\psi)\!-\!\frac{b_{c,lm}}{2\tau_{lm}^2}\cos(\psi)\right)\!V_l^2\\ &\qquad\quad\quad-
  \frac{Y_{lm}}{\tau_{lm}}V_l V_m\sin(\theta_{lm}-\delta_{lm}-\theta_{lm}^{shift}-\psi),\\
  \label{eq:general shifted OPF active power flow ml}
\nonumber &\!\!\!\tilde{P}_{ml} \!=\! \Re(\tilde{S}_{ml}) \!=\! \left(Y_{lm} \cos(\delta_{lm}+\psi)\!-\!\frac{b_{c,lm}}{2}  \sin(\psi)\right)\!V_m^2
  \\ &\qquad\quad\quad - \frac{Y_{lm}}{\tau_{lm}}V_m V_l \cos(\theta_{lm}+\delta_{lm}-\theta_{lm}^{shift}+\psi),\\
  \label{eq:general shifted OPF reactive power flow ml}
\nonumber &\!\!\!\tilde{Q}_{ml} \!=\! \Im(\tilde{S}_{ml}) \!=\!  \left(-Y_{lm}  \sin(\delta_{lm}+\psi)\!-\!\frac{b_{c,lm}}{2} \cos(\psi)\right)\!V_m^2 \\ &\qquad\quad\quad + \frac{Y_{lm}}{\tau_{lm}}V_m V_l \sin(\theta_{lm}+\delta_{lm}-\theta_{lm}^{shift}+\psi).
\end{align}
\end{subequations}

The  arguments of the trigonometric terms in~\eqref{eq:general rotated power flow} are not independent since $\cos(\theta_{lm}+\delta_{lm}-\theta_{lm}^{shift}+\psi)$ and $\sin(\theta_{lm}+\delta_{lm}-\theta_{lm}^{shift}+\psi)$ are linearly related with $\cos(\theta_{lm}-\delta_{lm}-\theta_{lm}^{shift}-\psi)$ and $\sin(\theta_{lm}-\delta_{lm}-\theta_{lm}^{shift}-\psi)$ via the general form of~\eqref{eq:terminal_relationships}. Extending~\eqref{eq:terminal_relationships} to consider off-nominal voltage ratios and non-zero phase shifts is accomplished by replacing $\theta_{lm}$ in~\eqref{eq:terminal_relationships} with $\theta_{lm}-\theta_{lm}^{shift}$.

Extensions of the expressions for the squared magnitudes of the current flows in the original QC relaxation~\eqref{eq:QC relaxation} and the RQC relaxation~\eqref{eq:RQC OPF}, $\ell_{lm}$ and $\tilde{\ell}_{lm}$, respectively, are derived by dividing $(P_{lm}^{2}+Q_{lm}^{2})$ and $(\tilde{P}_{lm}^{2}+\tilde{Q}_{lm}^{2})$ by $V_l^2$:
%
\begin{align}
\label{eq:squared current in QC relaxation}
\nonumber &\ell_{lm} = \left(\frac{Y_{lm}^2}{\tau_{lm}^4} - \frac{b_{c,lm}^2}{4\tau_{lm}^4} \right)V_l^2 + \frac{Y_{lm}^2}{\tau_{lm}^2}V_m^2 - \frac{b_{c,lm}}{\tau_{lm}^2}Q_{lm} \\
&\quad\quad - 2\frac{Y_{lm}^2}{\tau_{lm}^3}\left(\cos(\delta_{lm}) c_{lm} + \sin(\delta_{lm})s_{lm}\right),\\
%
%
%
%
%
%
%
%
%
&\nonumber \tilde{\ell}_{lm} = \bigg(\frac{Y_{lm}^2}{\tau_{lm}^4}+\frac{b_{c,lm}^2}{4\tau_{lm}^4}- \frac{Y_{lm}}{\tau_{lm}^4}b_{c,lm}\cos(\delta_{lm}+\psi)\sin(\psi)\\&\quad\quad\nonumber
 + \frac{Y_{lm}}{\tau_{lm}^4}b_{c,lm}\sin(\delta_{lm}+\psi)\cos(\psi)\bigg)V_l^2+\frac{Y_{lm}^2}{\tau_{lm}^2}V_m^2\\&\quad\quad\nonumber+\left(\frac{Y_{lm}}{\tau_{lm}^3}b_{c,lm}(\sin(\psi)- \frac{2Y_{lm}^2}{\tau_{lm}^3}\cos(\delta_{lm}+\psi)\right)\tilde{c}_{lm}\\&\quad\quad+\left(\frac{Y_{lm}}{\tau_{lm}^3}b_{c,lm}\cos(\psi)+\frac{2Y_{lm}^2}{\tau_{lm}^3}\sin(\delta_{lm}+\psi)\right)\tilde{s}_{lm}
\end{align}

Extending the TRQC relaxation~\eqref{eq:TRQC OPF} to the more general line model is achieved by modifying~\eqref{eq:angle_sumdiff3}:
\begingroup \small
\begin{align}
\label{eq:angle_sumdiff4}
& \begin{bmatrix}
\cos(\theta_{lm})\\
\sin(\theta_{lm}) 
\end{bmatrix} = M^{\prime}_{lm}
\begin{bmatrix} 
\sin(\theta_{lm} - \delta_{lm} - \psi-\theta_{lm}^{shift})\\ 
\cos(\theta_{lm} - \delta_{lm} - \psi-\theta_{lm}^{shift})
\end{bmatrix},
\end{align}\endgroup
where the constant matrix $M^{\prime}_{lm}$ is defined as
\begin{align*}
\footnotesize
M_{lm}^{\prime} = & \frac{1}{2}\left(\begin{bmatrix}
-\sin(\hat{\delta}_{lm}+\theta_{lm}^{shift}) & \cos(\hat{\delta}_{lm}+\theta_{lm}^{shift}) \\
\cos(\hat{\delta}_{lm}+\theta_{lm}^{shift}) & \sin(\hat{\delta}_{lm}+\theta_{lm}^{shift})
\end{bmatrix}\right. \\ &\!\!\!\!\!\!\!\!\!\!\!\!\!\!\!\!\!\!
\left. 
+\begin{bmatrix}
\sin(\hat{\delta}_{lm} -\theta_{lm}^{shift}) & \cos(\hat{\delta}_{lm}-\theta_{lm}^{shift}) \\
\cos(\hat{\delta}_{lm}-\theta_{lm}^{shift}) & -\sin(\hat{\delta}_{lm}-\theta_{lm}^{shift})
\end{bmatrix}\begin{bmatrix}
\alpha_{lm} & \beta_{lm} \\
-\beta_{lm} & \alpha_{lm}
\end{bmatrix} 
 \right)%
\end{align*}%
and, for notational convenience, $\hat{\delta}_{lm} = \delta_{lm} + \psi$.

\subsection{Parallel Lines}
\label{appendix:parallel}
In the original QC relaxation~\eqref{eq:QC relaxation}, the power flow equations for parallel lines between buses $l$ and $m$ shared the same envelopes, $\left\langle\cos(\theta_{lm})\right\rangle^C$ and $\left\langle\sin(\theta_{lm})\right\rangle^S$. In the RQC relaxation~\eqref{eq:RQC OPF}, the arguments of the trigonometric terms for parallel lines can differ due to the inclusion of the $\delta_{lm}$ terms. Rather than defining separate envelopes, we derive a linear relationship between the trigonometric terms for parallel lines. Let $\delta_{lm_1}$, $\delta_{lm_2}$ and $\theta_{lm_1}^{shift}$, $\theta_{lm_2}^{shift}$ be the admittance angles and phase shifts, respectively, for two parallel lines between buses~$l$ and~$m$. Applying the angle sum identity yields
\begin{equation}\small
\label{eq:angle_sum_parallel}
\begin{bmatrix}
\sin(\sigma_{lm_1} - \theta_{lm}) \\
\cos(\sigma_{lm_1} - \theta_{lm}) \\
\sin(\sigma_{lm_2} - \theta_{lm}) \\
\cos(\sigma_{lm_2} - \theta_{lm})
\end{bmatrix} = \begin{bmatrix}
\sin(\sigma_{lm_1}) & -\cos(\sigma_{lm_1}) \\
\cos(\sigma_{lm_1}) & \hphantom{-}\sin(\sigma_{lm_1}) \\
\sin(\sigma_{lm_2}) & -\cos(\sigma_{lm_2}) \\
\cos(\sigma_{lm_2}) & \hphantom{-}\sin(\sigma_{lm_2})
\end{bmatrix}\begin{bmatrix}
\cos(\theta_{lm}) \\
\sin(\theta_{lm})
\end{bmatrix},
\end{equation}
where, for notational convenience, $\sigma_{lm_1} = \delta_{lm_1} + \theta_{lm_1}^{shift} + \psi$ and $\sigma_{lm_2} = \delta_{lm_2} + \theta_{lm_2}^{shift} + \psi$. Rearranging~\eqref{eq:angle_sum_parallel} to eliminate $\cos(\theta_{lm})$ and $\sin(\theta_{lm})$ yields the desired linear relationship: 
\begingroup \small
\begin{align}
\nonumber
& \begin{bmatrix}
\sin(\theta_{lm} - \sigma_{lm_2}) \\
\cos(\theta_{lm} - \sigma_{lm_2})
\end{bmatrix}  = \\ \label{eq:terminal_relationships_parallel} & \; \begin{bmatrix}
\hphantom{-}\cos(\sigma_{lm_1}-\sigma_{lm_2}) & \sin(\sigma_{lm_1}-\sigma_{lm_2}) \\
-\sin(\sigma_{lm_1}-\sigma_{lm_2}) & \cos(\sigma_{lm_1}-\sigma_{lm_2})
\end{bmatrix}\begin{bmatrix}
\sin(\theta_{lm} - \sigma_{lm_1})\\
\cos(\theta_{lm} - \sigma_{lm_1})
\end{bmatrix}.
\end{align}\endgroup%
Since the matrix in~\eqref{eq:terminal_relationships_parallel} is invertible, this relationship is always well defined.

\subsection{Tighter Boundaries for Certain Trigonometric Envelopes}
\label{appendix:Proof of Tighter Lower Bound of RQC relaxation}
This appendix formalizes and proves a statement in Section~\ref{Different envelope epresentation explanations} regarding the tightness of the trigonometric envelopes in the original formulation of the QC relaxation~\eqref{eq:QC relaxation} and the proposed polar admittance QC formulation~\eqref{eq:RQC OPF}.

To assist the derivations in this appendix, we define a function $F(\theta_{lm})$ which represents the difference between the trigonometric function $\cos(\theta_{lm} - \delta_{lm})$ itself and the line which connects the endpoints of $\cos(\theta_{lm} - \delta_{lm})$ at $\theta_{lm}^{min}$ and $\theta_{lm}^{max}$:
\begin{align}
\label{eq:F}
\nonumber \noindent F(\theta_{lm}) & = \cos(\theta_{lm} - \delta_{lm}) - \cos(\theta_{lm}^{max}-\delta_{lm})\\
&-\frac{\cos(\theta_{lm}^{max}-\delta_{lm})-\cos(\theta_{lm}^{min}-\delta_{lm})}{\theta_{lm}^{max}-\theta_{lm}^{min}}\left(\theta_{lm}-\theta_{lm}^{max}\right)
\end{align}
Fig.~\ref{fig:proof} shows illustrative examples of the function $Y_{lm}\cos(\theta_{lm} - \delta_{lm})$ (black curve) and the line connecting the endpoints of this function at $\theta_{lm}^{min}$ and $\theta_{lm}^{max}$ (dashed red line) on the left, with corresponding visualizations of the function $F(\theta_{lm})$ itself on the right.

\begin{figure}[ht]
\def\tabularxcolumn#1{m{#1}}
\subfloat[$-165^\circ \leq \theta_{lm}-\delta_{lm} \leq 15^\circ$\label{fig:orig1}]{
\begin{tabular}{cc}
\small $\cos(\theta_{lm} - \delta_{lm})$ \vspace*{-0.1em} & \small $F(\theta_{lm})$ \vspace*{-0.1em}\\
\includegraphics[width=4.1cm]{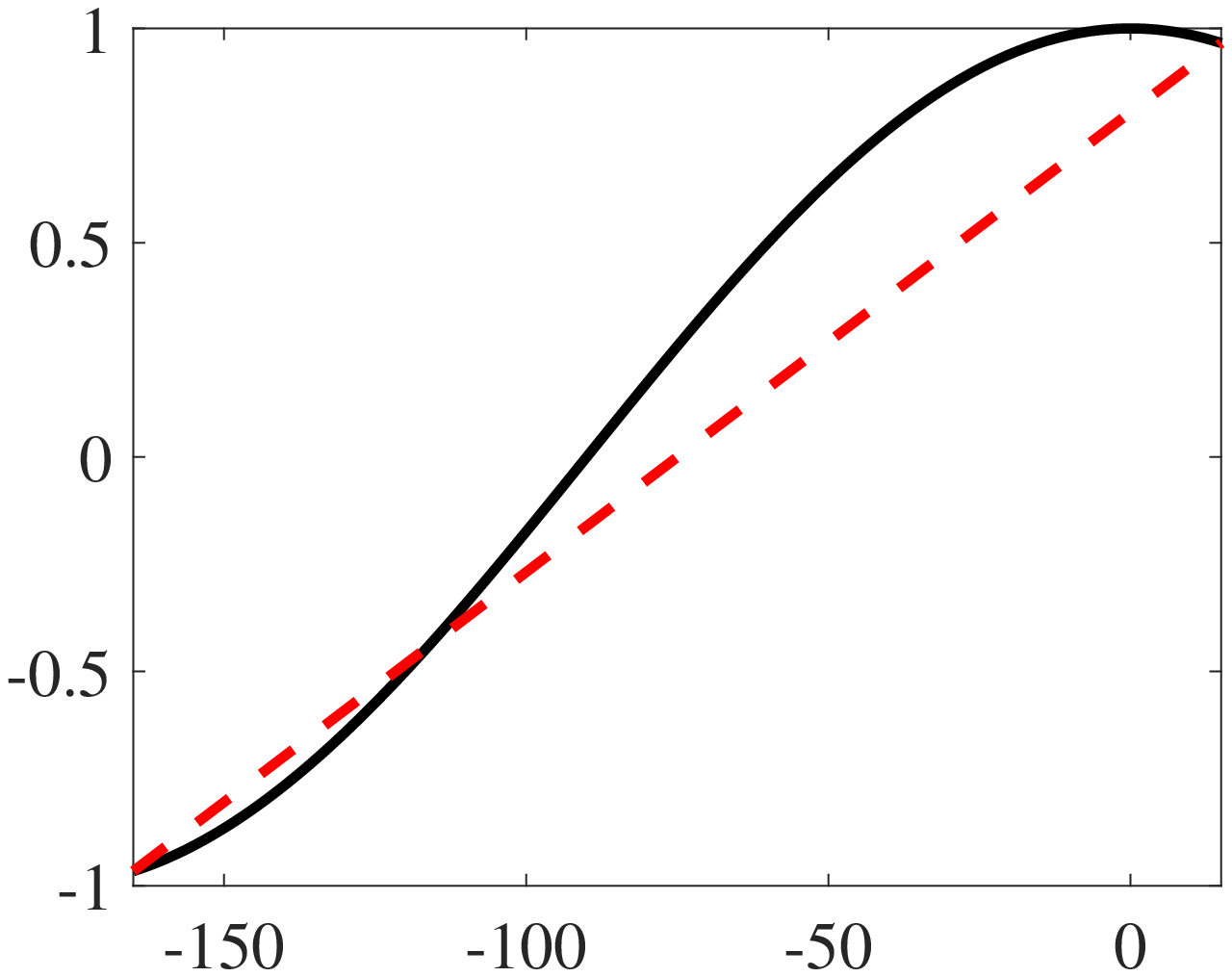} & \includegraphics[width=4.1cm]{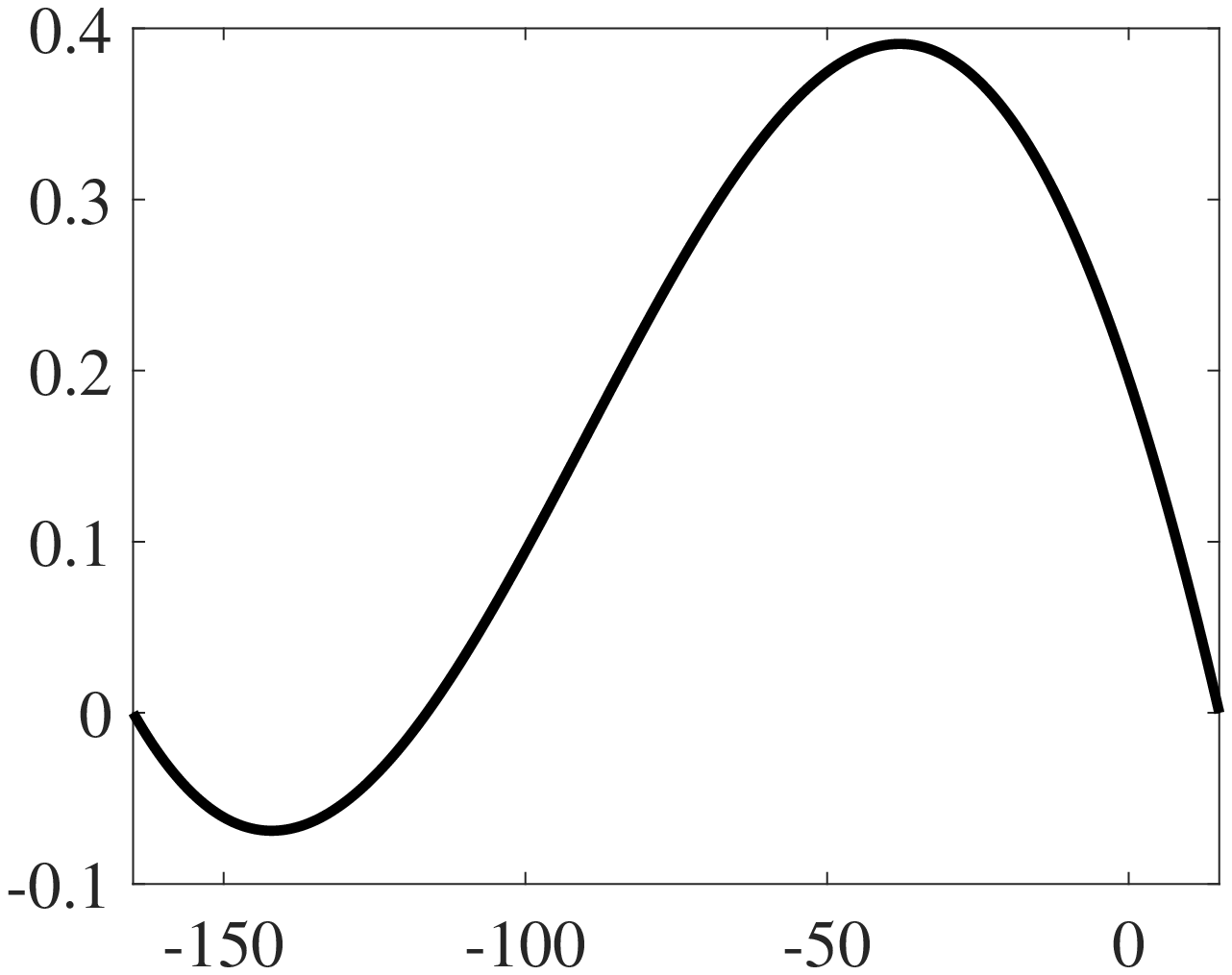}
\end{tabular}%
}\\
\subfloat[$-165^\circ \leq \theta_{lm}-\delta_{lm} \leq 15^\circ$\label{fig:orig3}]{
\begin{tabular}{cc}
\small $\cos(\theta_{lm} - \delta_{lm})$ \vspace*{-0.1em} & \small $F(\theta_{lm})$ \vspace*{-0.1em} \\
\includegraphics[width=4.1cm]{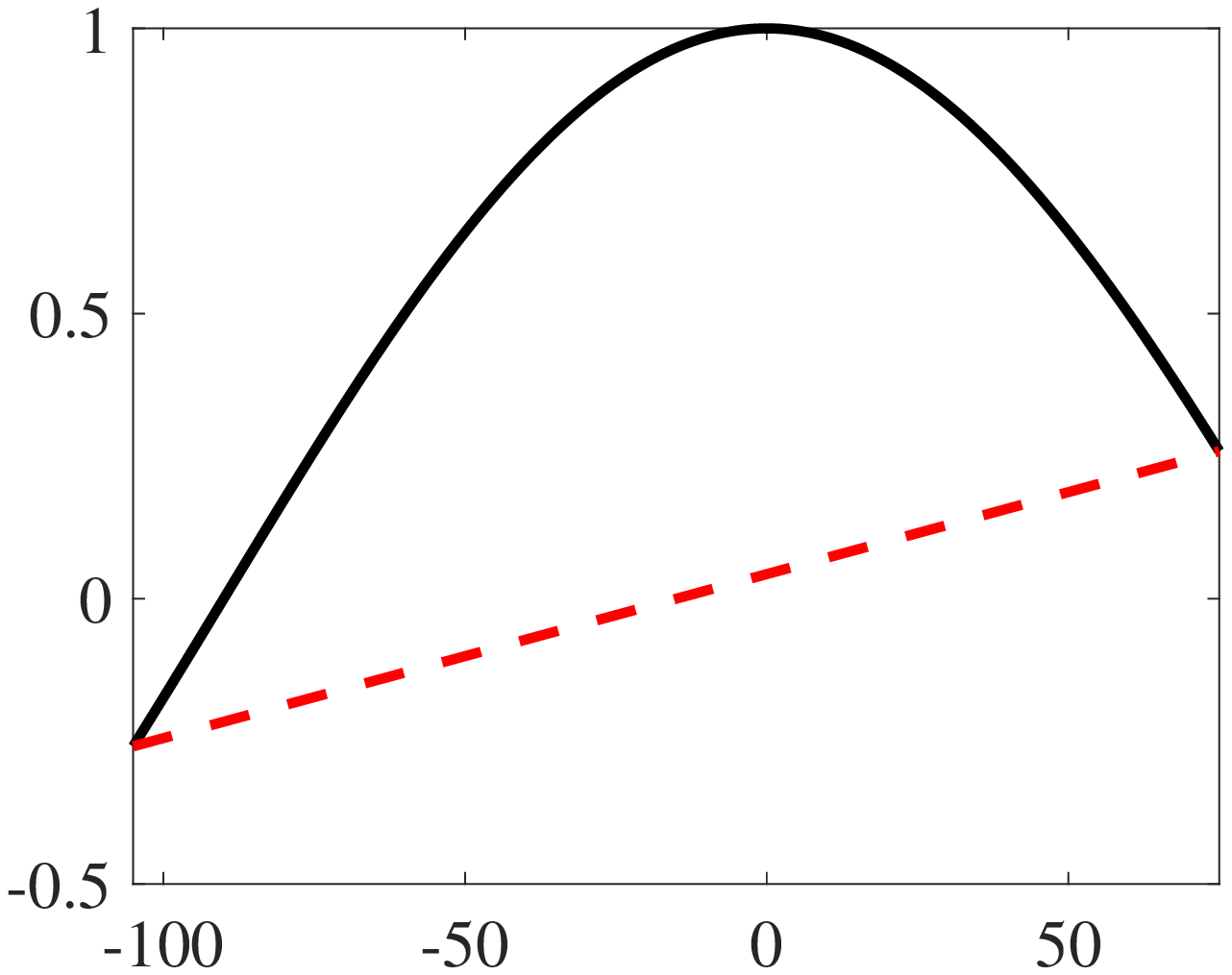} & \includegraphics[width=4.1cm]{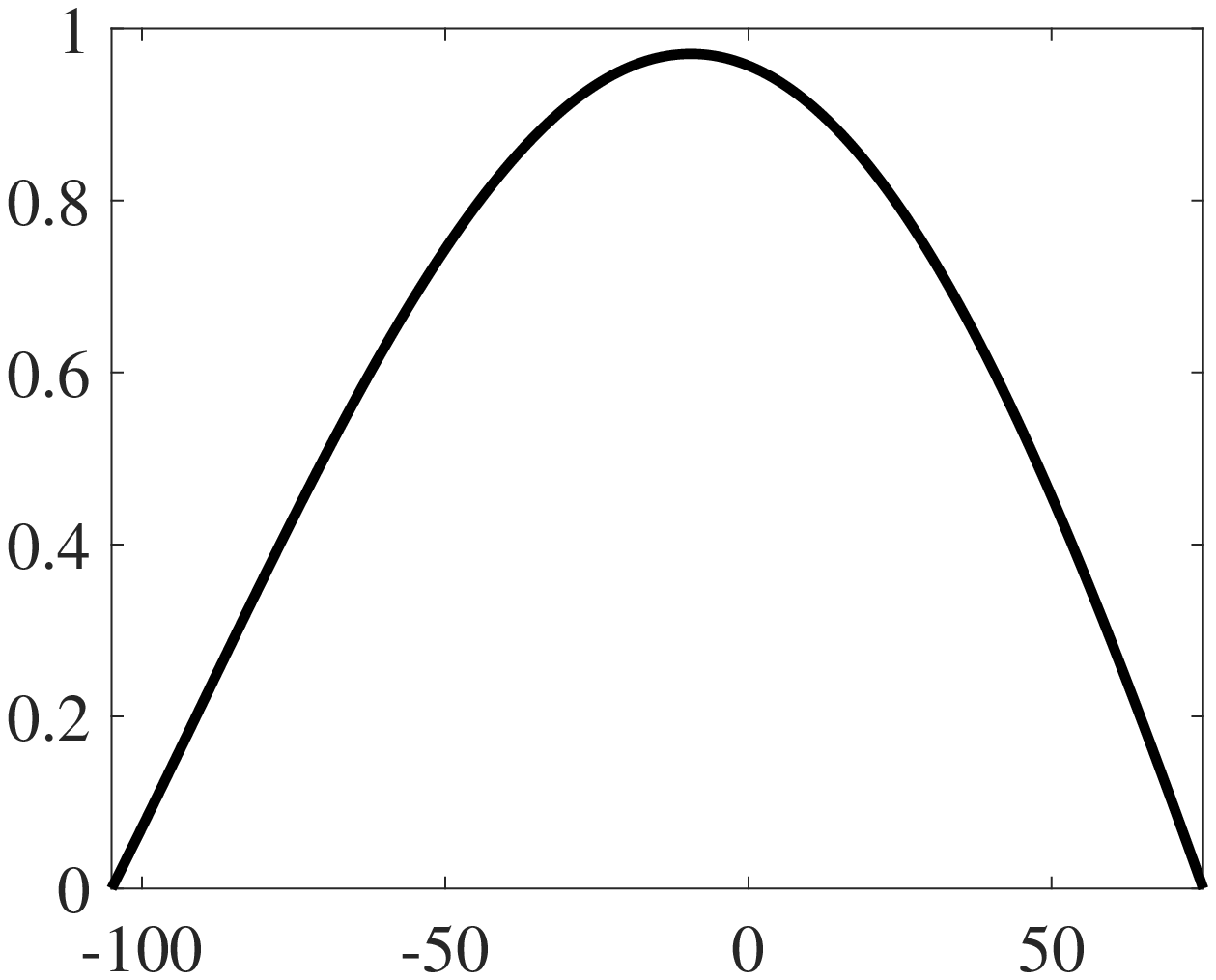}
\end{tabular}%
}
\centering
\caption{The left figures show visualizations of the function $\cos(\theta_{lm} - \delta_{lm})$ (black curve) and the line connecting the endpoints of this function at $\theta_{lm}^{min}$ and $\theta_{lm}^{max}$ (dashed red line) for different values of $\delta_{lm}$, $\theta_{lm}^{min}$, and $\theta_{lm}^{max}$. The right  figures show the corresponding function $F(\theta_{lm})$.}
\label{fig:proof}  \vspace{0.5cm}
\end{figure}

The derivative of $F(\theta_{lm})$ is
\begin{align}
\nonumber & \frac{dF(\theta_{lm})}{d \theta_{lm}} = -\sin(\theta_{lm}-\delta_{lm}) \\ 
& \quad - \frac{\cos(\theta_{lm}^{max} -\delta_{lm})-\cos(\theta_{lm}^{min}-\delta_{lm})}{\theta_{lm}^{max}-\theta_{lm}^{min}}.
\end{align}

A key quantity in the following proposition is the set of zeros of the derivative of $F(\theta_{lm})$, i.e., the set of solutions to $\frac{dF(\theta_{lm})}{d \theta_{lm}} = 0$. This set, which we denote by $\mathcal{Z}_{\theta_{lm}^{min}, \theta_{lm}^{max}, \delta_{lm}}$ where the subscripts indicate that the set is parameterized by $\theta_{lm}^{min}$, $\theta_{lm}^{max}$, and $\delta_{lm}$, is
\begin{align}
\nonumber & \mathcal{Z}_{\theta_{lm}^{min}, \theta_{lm}^{max}, \delta_{lm}} = \\ 
\nonumber & \Big\lbrace (-1)^\kappa \arcsin\left(\frac{\cos(\theta_{lm}^{min}-\delta_{lm}) - \cos(\theta_{lm}^{max}-\delta_{lm})}{\left(\theta_{lm}^{max}-\theta_{lm}^{min}\right)}\right) + \pi\kappa,\\
\nonumber 
& \qquad \kappa = \ldots, -3, -2, -1, 0, 1, 2, 3, \ldots \Big\rbrace.
\end{align}
Finally, we let $\left|\,\cdot\,\right|$ denote the cardinality of a set.

Using these definitions, we next state and prove the following proposition.
\begin{proposition}
The lower boundary of the envelope $Y_{lm}\left\langle \cos(\theta_{lm} - \delta_{lm}) \right\rangle^C$ is at least as tight as the lower boundary of the envelope $g_{lm}\left\langle\cos(\theta_{lm})\right\rangle^C + b_{lm}\left\langle\sin(\theta_{lm})\right\rangle^S$ if $\theta_{lm}^{min}$, $\theta_{lm}^{max}$, and $\delta_{lm}$ satisfy both of the following conditions:
\begin{subequations}
\label{eq:lower_bound_condition}
\begin{align}
\label{eq:lower_bound_condition_1}
& \left|\mathcal{Z}_{\theta_{lm}^{min}, \theta_{lm}^{max}, \delta_{lm}} \bigcap \left\lbrace \theta_{lm}^{min} < \theta_{lm} < \theta_{lm}^{max}\right\rbrace \right| = 1,\\
\label{eq:lower_bound_condition_2}
& F\left(\left(\theta_{lm}^{max} + \theta_{lm}^{min}\right) / 2\right) > 0.
\end{align}
\end{subequations}
Moreover, the upper boundary of the envelope $Y_{lm}\left\langle \cos(\theta_{lm} - \delta_{lm}) \right\rangle^C$ is at least as tight as the upper boundary of the envelope $g_{lm}\left\langle\cos(\theta_{lm})\right\rangle^C + b_{lm}\left\langle\sin(\theta_{lm})\right\rangle^S$ if $\theta_{lm}^{min}$, $\theta_{lm}^{max}$, and $\delta_{lm}$ satisfy both~\eqref{eq:lower_bound_condition_1} and the condition
\begin{equation}
\label{eq:upper_bound_condition_2}
F\left(\left(\theta_{lm}^{max} + \theta_{lm}^{min}\right) / 2\right) < 0.
\end{equation}

\begin{proof}
The proof is based on the following observation: if the line connecting the points $(\theta_{lm}^{min},\,\cos(\theta_{lm}^{min} - \delta_{lm}))$ and $(\theta_{lm}^{max},\,\cos(\theta_{lm}^{max} - \delta_{lm}))$ (i.e., the dashed red line in Fig.~\ref{fig:proof}) does not intersect the function $\cos(\theta_{lm} - \delta_{lm})$ itself within the range $\theta_{lm}^{min} < \theta_{lm} < \theta_{lm}^{max}$, then this line is either the lower boundary or upper boundary of the tighest convex envelope for the function $\cos(\theta_{lm} - \delta_{lm})$ within this range. (For instance, the dashed red line in Fig.~\ref{fig:orig3} is the lower boundary of the tighest envelope for $\cos(\theta_{lm} - \delta_{lm})$ within the range $-165^\circ \leq \theta_{lm} \leq 15^\circ$.) In this case, the line is the tightest lower (upper) boundary if the function $\cos(\theta_{lm} - \delta_{lm})$ is above (below) the line for any point between $\theta_{lm}^{min}$ and $\theta_{lm}^{max}$ (e.g., the midpoint $(\theta_{lm}^{min} + \theta_{lm}^{max})/2$, which is used in~\eqref{eq:lower_bound_condition_2} and~\eqref{eq:upper_bound_condition_2}).

Observe that the line connecting the points $(\theta_{lm}^{min},\,\cos(\theta_{lm}^{min} - \delta_{lm}))$ and $(\theta_{lm}^{max},\,\cos(\theta_{lm}^{max} - \delta_{lm}))$ does not intersect the function $\cos(\theta_{lm} - \delta_{lm})$ between $\theta_{lm}^{min}$ and $\theta_{lm}^{max}$ if and only if $F(\theta_{lm})$ is non-zero for all $\theta_{lm}^{min} < \theta_{lm} < \theta_{lm}^{max}$. We next argue that this is implied by~\eqref{eq:lower_bound_condition_1}.

The condition~\eqref{eq:lower_bound_condition_1} is equivalent to the existence of one critical point $\theta_{lm}^\ast$ of the function $F(\theta_{lm})$. (i.e., the derivative of $F(\theta)$ has a single zero, $\theta_{lm}^\ast$, in the range $\theta_{lm}^{min} < \theta_{lm} < \theta_{lm}^{max}$. Since $F(\theta_{lm})$ is continuous and $F(\theta_{lm}^{min}) = F(\theta_{lm}^{max}) = 0$, the critical point $\theta_{lm}^\ast$ must either correspond to a minimum or maximum of $F(\theta_{lm})$. Since the function $F(\theta_{lm})$ is zero at the endpoints $\theta_{lm}^{min}$ and $\theta_{lm}^{max}$,  having a single minimum or maximum in the range $\theta_{lm}^{min} < \theta_{lm} < \theta_{lm}^{max}$ implies that $F(\theta_{lm}) \neq 0$ within this range. 



To complete the conditions in the proposition,~\eqref{eq:lower_bound_condition_2} and~\eqref{eq:upper_bound_condition_2} determine whether the line connecting the points $(\theta_{lm}^{min},\,\cos(\theta_{lm}^{min} - \delta_{lm}))$ and $(\theta_{lm}^{max},\,\cos(\theta_{lm}^{max} - \delta_{lm}))$ is above or below the function $\cos(\theta_{lm} - \delta_{lm})$ by evaluating the function $F(\theta_{lm})$ at an arbitrary point between $\theta_{lm}^{min}$ and $\theta_{lm}^{max}$, here selected to be the midpoint $\left(\theta_{lm}^{min} + \theta_{lm}^{max}\right)/2$.


Since multiplication by $Y_{lm}$ only rescales (but does not otherwise change) the envelope $\left\langle \cos(\theta_{lm} - \delta_{lm})\right\rangle^C$, the arguments above trivially extend to $Y_{lm}\left\langle \cos(\theta_{lm} - \delta_{lm})\right\rangle^C$. Moreover, since $Y_{lm}\cos(\theta_{lm} - \delta_{lm}) = g_{lm}\cos(\theta_{lm}) + b_{lm}\sin(\theta_{lm})$, the envelope $g_{lm}\left\langle\cos(\theta_{lm})\right\rangle^C + b_{lm}\left\langle \sin(\theta_{lm})\right\rangle^S$ cannot be tighter than the tightest possible envelope for $Y_{lm}\cos(\theta_{lm} - \delta_{lm})$. Since the boundaries of $Y_{lm}\left\langle \cos(\theta_{lm} - \delta_{lm})\right\rangle^C$ considered in the proof form portions of the tightest possible convex envelope for $Y_{lm}\cos(\theta_{lm} - \delta_{lm})$, they are at least as tight as the corresponding boundaries of the envelope $g_{lm}\left\langle\cos(\theta_{lm})\right\rangle^C + b_{lm}\left\langle \sin(\theta_{lm})\right\rangle^S$. Furthermore, the example envelopes in Fig.~\ref{fig:psi_impact} show that the corresponding boundaries of $Y_{lm}\left\langle \cos(\theta_{lm} - \delta_{lm})\right\rangle^C$ are strictly tighter than those of $g_{lm}\left\langle\cos(\theta_{lm})\right\rangle^C + b_{lm}\left\langle \sin(\theta_{lm})\right\rangle^S$ for some values of $\delta_{lm}$, $\theta_{lm}^{min}$, and $\theta_{lm}^{max}$.

\end{proof}
\end{proposition}


We finally note that values of $\theta_{lm}^{min}$, $\theta_{lm}^{max}$, and $\delta_{lm}$ such that $\max(-90^\circ, -90^\circ + \delta_{lm}) \leq\theta_{lm}^{min} < \theta_{lm}^{max} \leq \min(90^\circ, 90^\circ + \delta_{lm})$ satisfy~\eqref{eq:lower_bound_condition}. Thus, the trigonometric envelopes corresponding to the polar admittance representation have lower boundaries that are at least as tight as those in the original QC relaxation for many typical values of $\theta_{lm}^{min}$, $\theta_{lm}^{max}$, and~$\delta_{lm}$.

\subsection{Expressions for the Vertices of the Polytope Associated with the Trilinear Products}
\label{apx:intersection}

This appendix presents expressions for the vertices of the polytope consisting of the black dashed lines in Fig.~\ref{fig_intersection}.
To compute the coordinates of these vertices (black dots in Fig.~\ref{fig_intersection_points}), we intersect the edges of the receiving end polytope, which is formed by the upper and lower bounds on the receiving end quantities, $\vphantom{\tilde{S}_{lm}}\smash{\overline{\tilde{S}}}_{lm}^{(r)}$, $\vphantom{\tilde{C}_{lm}}\smash{\overline{\tilde{C}}}_{lm}^{(r)}$ and $\vphantom{\tilde{S}_{lm}}\smash{\underline{\tilde{S}}}_{lm}^{(r)}$, $\vphantom{\tilde{C}_{lm}}\smash{\underline{\tilde{C}}}_{lm}^{(r)}$, respectively, with the edges of the sending end polytope, which is formed by the upper and lower bounds on the sending end quantities $\vphantom{\tilde{S}_{lm}}\smash{\overline{\tilde{S}}}_{lm}^{(s)}$, $\vphantom{\tilde{C}_{lm}}\smash{\overline{\tilde{C}}}_{lm}^{(s)}$ and $\vphantom{\tilde{S}_{lm}}\smash{\underline{\tilde{S}}}_{lm}^{(s)}$, $\vphantom{\tilde{C}_{lm}}\smash{\underline{\tilde{C}}}_{lm}^{(s)}$, respectively. 

When written in terms of the sending end quantities $\tilde{S}_{lm}^{(s)}$ and $\tilde{C}_{lm}^{(s)}$, the coordinates for the upper and lower bounds on the receiving end quantities are functions of $\psi$. To write the coordinates of the vertices as functions of $\psi$, consider the line segments labeled in Fig.~\ref{fig_intersection_points}. The yellow and purple polytopes in this figure represent the bounds on the sending and receiving end quantities, respectively. Table~\ref{tab:line_segments_table} describes the relevant intersections of the line segments that form these polytopes. For the ranges of $\psi$ in the first column of Table~\ref{tab:line_segments_table}, the remaining columns indicate the line segments whose intersections form the corresponding vertices. The coordinates of these intersections are given in Table~\ref{tab:intersection}. As an example for $-45^\circ \le\psi\le0^\circ$, the $A^{\prime}D^{\prime}$ line segment in Fig.~\ref{fig_intersection_points} should intersect line segments AB and AD. The coordinates of these intersections are given in rows 13 and 16 of Table~\ref{tab:intersection}.

\begin{table}
\caption{Line segment intersections corresponding to Fig.~\ref{fig_intersection_points}}
\footnotesize
\label{tab:line_segments_table}
\begin{tabular}{|c|c|c|c|c|}
\hline 
$\psi$ (degrees) & $A^{\prime}B^{\prime}$ & $B^{\prime}C^{\prime}$ & $C^{\prime}D^{\prime}$  & $A^{\prime}D^{\prime}$
 \tabularnewline
\hline 
\hline 
\!\!\!\!\!${-}45\le\psi\le 0$\!\!\!\!\! & AB \& BC & BC \& CD & CD \& AD & AB \& AD \tabularnewline
\hline 
 \!\!\!$\!\!{-}90\le\psi\le {-}45$\!\!\!& BC \& CD & CD \& AD & AB \& AD  & AB \& BC\tabularnewline
\hline 
 \!\!\!${-}135\le\psi\le {-}90$\!\!\!& CD \& AD & AB \& AD  &AB \& BC  & BC \& CD\tabularnewline
\hline 
\!\!\!\!\! ${-}180\le\psi\le {-}135$\!\!\!\!& AB \& AD & AB \& BC & BC \& CD & CD \& AD\tabularnewline
\hline 
 \!\!\!\!\!$0\le\psi\le45$\!\!\!\!\!&  AD \& AB& AB \& BC & BC \& CD & CD \& AD \tabularnewline
\hline 
 \!\!\!\!\!$45\le\psi\le90$\!\!\!\!\!&  CD \& AD& AD \& AB & AB \& BC & BC \& CD \tabularnewline
\hline 
 \!\!\!\!\!$90\le\psi\le135$\!\!& BC \& CD & CD \& AD & AD \& AB & AB \& BC \tabularnewline
\hline 
 \!\!\!\!\!$135\le\psi\le180$\!\!\!\!\!& AB \& BC  & BC \& CD & CD \& AD & AD \& AB\tabularnewline
\hline 
\end{tabular}
\end{table}

\begin{table}
\caption{Coordinates of the line segment intersections in Table~\ref{tab:line_segments_table}.}
\begin{small}
\label{tab:intersection}
\begin{tabular}{|c|c|}
\hline 
 \!\!\!Line Segments\!\!\!& Coordinates of the Intersection Point \tabularnewline
\hline 
\hline 
 $A^{\prime}B^{\prime}$ \& AB& $\left(\frac{\alpha_{lm}}{\beta_{lm}}\vphantom{\tilde{S}_{lm}}\smash{\overline{\tilde{S}}}_{lm}^{(s)}-\beta_{lm}\vphantom{\tilde{S}_{lm}}\smash{\overline{\tilde{S}}}_{lm}^{(r)})-\frac{\alpha_{lm}^2\vphantom{\tilde{S}_{lm}}\smash{\overline{\tilde{S}}}_{lm}^{(r)})}{\beta_{lm}},\;\vphantom{\tilde{S}_{lm}}\smash{\overline{\tilde{S}}}_{lm}^{(s)}\right)$ \tabularnewline
\hline 
  $A^{\prime}B^{\prime}$ \& BC& $\left(\vphantom{\tilde{C}_{lm}}\smash{\underline{\tilde{C}}}_{lm}^{(s)},\; \frac{\beta_{lm}\vphantom{\underline{\tilde{C}}_{lm}}\smash{\underline{\tilde{C}}}_{lm}^{(s)})}{\alpha_{lm}}+\frac{\beta_{lm}^2\vphantom{\tilde{S}_{lm}}\smash{\overline{\tilde{S}}}_{lm}^{(r)})}{\alpha_{lm}}+ \alpha_{lm}\vphantom{\tilde{S}_{lm}}\smash{\overline{\tilde{S}}}_{lm}^{(r)}\right)$ \tabularnewline
\hline 
  $A^{\prime}B^{\prime}$ \& CD& $\left(\frac{\alpha_{lm}}{\beta_{lm}}\vphantom{\tilde{S}_{lm}}\smash{\underline{\tilde{S}}}_{lm}^{(s)}-\beta_{lm}\vphantom{\tilde{S}_{lm}}\smash{\overline{\tilde{S}}}_{lm}^{(r)})-\frac{\alpha_{lm}^2\vphantom{\tilde{S}_{lm}}\smash{\overline{\tilde{S}}}_{lm}^{(r)})}{\beta_{lm}},\;\underline{\tilde{S}}_{lm}^{(s)}\right)$\tabularnewline
\hline 
  $A^{\prime}B^{\prime}$ \& AD& $\left(\vphantom{\tilde{C}_{lm}}\smash{\overline{\tilde{C}}}_{lm}^{(s)},\;\frac{\beta_{lm}\vphantom{\tilde{C}_{lm}}\smash{\overline{\tilde{C}}}_{lm}^{(s)})}{\alpha_{lm}}+\frac{\beta_{lm}^2\vphantom{\tilde{S}_{lm}}\smash{\overline{\tilde{S}}}_{lm}^{(r)})}{\alpha_{lm}}+ \alpha_{lm}\vphantom{\tilde{S}_{lm}}\smash{\overline{\tilde{S}}}_{lm}^{(r)}\right)$\tabularnewline
\hline 
 $B^{\prime}C^{\prime}$ \& BC& $\left(\vphantom{\tilde{C}_{lm}}\smash{\underline{\tilde{C}}}_{lm}^{(s)},\;-\frac{\alpha_{lm}\vphantom{\tilde{C}_{lm}}\smash{\underline{\tilde{C}}}_{lm}^{(s)})}{\beta_{lm}}+\frac{\alpha_{lm}^2\vphantom{\tilde{C}_{lm}}\smash{\underline{\tilde{C}}}_{lm}^{(r)})}{\beta_{lm}}+ \beta_{lm}\vphantom{\tilde{C}_{lm}}\smash{\underline{\tilde{C}}}_{lm}^{(r)}\right)$\tabularnewline
\hline 
 $B^{\prime}C^{\prime}$ \& CD& $\left(-\frac{\beta_{lm}}{\alpha_{lm}}\vphantom{\tilde{S}_{lm}}\smash{\underline{\tilde{S}}}_{lm}^{(s)}+\alpha_{lm}\vphantom{\tilde{C}_{lm}}\smash{\underline{\tilde{C}}}_{lm}^{(r)})+\frac{\beta_{lm}^2\vphantom{\tilde{C}_{lm}}\smash{\underline{\tilde{C}}}_{lm}^{(r)})}{\alpha_{lm}},\;\vphantom{\tilde{S}_{lm}}\smash{\underline{\tilde{S}}}_{lm}^{(s)}\right)$\tabularnewline
\hline 
 $B^{\prime}C^{\prime}$ \& AB& $\left(-\frac{\beta_{lm}}{\alpha_{lm}}\vphantom{\tilde{S}_{lm}}\smash{\overline{\tilde{S}}}_{lm}^{(s)}+\alpha_{lm}\vphantom{\tilde{C}_{lm}}\smash{\underline{\tilde{C}}}_{lm}^{(r)})+\frac{\beta_{lm}^2\vphantom{\tilde{C}_{lm}}\smash{\underline{\tilde{C}}}_{lm}^{(r)})}{\alpha_{lm}},\; \vphantom{\tilde{S}_{lm}}\smash{\overline{\tilde{S}}}_{lm}^{(s)}\right)$\tabularnewline
\hline 
 $B^{\prime}C^{\prime}$ \& AD & $\left(\vphantom{\tilde{C}_{lm}}\smash{\overline{\tilde{C}}}_{lm}^{(s)},\;-\frac{\alpha_{lm}\vphantom{\tilde{C}_{lm}}\smash{\overline{\tilde{C}}}_{lm}^{(s)})}{\beta_{lm}}+\frac{\alpha_{lm}^2\vphantom{\tilde{C}_{lm}}\smash{\underline{\tilde{C}}}_{lm}^{(r)})}{\beta_{lm}}+ \beta_{lm}\vphantom{\tilde{C}_{lm}}\smash{\underline{\tilde{C}}}_{lm}^{(r)}\right)$ \tabularnewline
\hline 
  $C^{\prime}D^{\prime}$ \& AB& $\left(\frac{\alpha_{lm}}{\beta_{lm}}\vphantom{\tilde{S}_{lm}}\smash{\overline{\tilde{S}}}_{lm}^{(s)}-\beta_{lm}\vphantom{\tilde{S}_{lm}}\smash{\underline{\tilde{S}}}_{lm}^{(r)})-\frac{\alpha_{lm}^2\vphantom{\tilde{S}_{lm}}\smash{\underline{\tilde{S}}}_{lm}^{(r)})}{\beta_{lm}},\;\vphantom{\tilde{C}_{lm}}\smash{\overline{\tilde{S}}}_{lm}^{(s)}\right)$\tabularnewline
\hline 
   $C^{\prime}D^{\prime}$ \& BC& $\left(\vphantom{\tilde{C}_{lm}}\smash{\underline{\tilde{C}}}_{lm}^{(s)},\;\frac{\beta_{lm}\vphantom{\tilde{C}_{lm}}\smash{\underline{\tilde{C}}}_{lm}^{(s)})}{\alpha_{lm}}+\frac{\beta_{lm}^2\vphantom{\tilde{S}_{lm}}\smash{\underline{\tilde{S}}}_{lm}^{(r)})}{\alpha_{lm}}+ \alpha_{lm}\vphantom{\tilde{S}_{lm}}\smash{\underline{\tilde{S}}}_{lm}^{(r)}\right)$\tabularnewline
\hline 
   $C^{\prime}D^{\prime}$ \& CD& $\left(\frac{\alpha_{lm}}{\beta_{lm}}\vphantom{\tilde{C}_{lm}}\smash{\underline{\tilde{S}}}_{lm}^{(s)}-\beta_{lm}\vphantom{\tilde{S}_{lm}}\smash{\underline{\tilde{S}}}_{lm}^{(r)})-\frac{\alpha_{lm}^2\vphantom{\tilde{S}_{lm}}\smash{\underline{\tilde{S}}}_{lm}^{(r)})}{\beta_{lm}},\; \vphantom{\tilde{C}_{lm}}\smash{\underline{\tilde{S}}}_{lm}^{(s)}\right)$\tabularnewline
\hline 
   $C^{\prime}D^{\prime}$ \& AD& $\left(\vphantom{\tilde{C}_{lm}}\smash{\overline{\tilde{C}}}_{lm}^{(s)},\; \frac{\beta_{lm}\vphantom{\tilde{C}_{lm}}\smash{\overline{\tilde{C}}}_{lm}^{(s)})}{\alpha_{lm}}+\frac{\beta_{lm}^2\vphantom{\tilde{S}_{lm}}\smash{\underline{\tilde{S}}}_{lm}^{(r)})}{\alpha_{lm}}+ \alpha_{lm}\vphantom{\tilde{S}_{lm}}\smash{\underline{\tilde{S}}}_{lm}^{(r)}\right)$\tabularnewline
\hline 
   $A^{\prime}D^{\prime}$ \& AB& $\left(-\frac{\beta_{lm}}{\alpha_{lm}}\vphantom{\tilde{S}_{lm}}\smash{\overline{\tilde{S}}}_{lm}^{(s)}+\alpha_{lm}\vphantom{\tilde{C}_{lm}}\smash{\overline{\tilde{C}}}_{lm}^{(r)})+\frac{\beta_{lm}^2\vphantom{\tilde{C}_{lm}}\smash{\overline{\tilde{C}}}_{lm}^{(r)})}{\alpha_{lm}},\;\vphantom{\tilde{S}_{lm}}\smash{\overline{\tilde{S}}}_{lm}^{(s)}\right)$ \tabularnewline
\hline 
    $A^{\prime}D^{\prime}$ \& BC& $\left(\vphantom{\tilde{C}_{lm}}\smash{\underline{\tilde{C}}}_{lm}^{(s)},\;-\frac{\alpha_{lm}\vphantom{\tilde{C}_{lm}}\smash{\underline{\tilde{C}}}_{lm}^{(s)})}{\beta_{lm}}+\frac{\alpha_{lm}^2\vphantom{\tilde{C}_{lm}}\smash{\overline{\tilde{C}}}_{lm}^{(r)})}{\beta_{lm}}+ \beta_{lm}\vphantom{\tilde{C}_{lm}}\smash{\overline{\tilde{C}}}_{lm}^{(r)}\right)$ \tabularnewline
\hline 
    $A^{\prime}D^{\prime}$ \& CD& $\left(-\frac{\beta_{lm}}{\alpha_{lm}}\vphantom{\tilde{S}_{lm}}\smash{\underline{\tilde{S}}}_{lm}^{(s)}+\alpha_{lm}\vphantom{\tilde{C}_{lm}}\smash{\overline{\tilde{C}}}_{lm}^{(r)})+\frac{\beta_{lm}^2\vphantom{\tilde{C}_{lm}}\smash{\overline{\tilde{C}}}_{lm}^{(r)})}{\alpha_{lm}},\; \vphantom{\tilde{S}_{lm}}\smash{\underline{\tilde{S}}}_{lm}^{(s)}\right)$ \tabularnewline
\hline 
    $A^{\prime}D^{\prime}$ \& AD& $\left(\vphantom{\tilde{C}_{lm}}\smash{\overline{\tilde{C}}}_{lm}^{(s)},\; -\frac{\alpha_{lm}\vphantom{\tilde{C}_{lm}}\smash{\overline{\tilde{C}}}_{lm}^{(s)})}{\beta_{lm}}+\frac{\alpha_{lm}^2\vphantom{\tilde{C}_{lm}}\smash{\overline{\tilde{C}}}_{lm}^{(r)})}{\beta_{lm}}+ \beta_{lm}\vphantom{\tilde{C}_{lm}}\smash{\overline{\tilde{C}}}_{lm}^{(r)}\right)$  \tabularnewline
\hline 
\end{tabular}
\end{small}
\vspace*{-1em}
\end{table}

\begin{figure}[t!]
    \centering
\begin{tikzpicture}
\node at (0,0) {\includegraphics[scale=0.37]{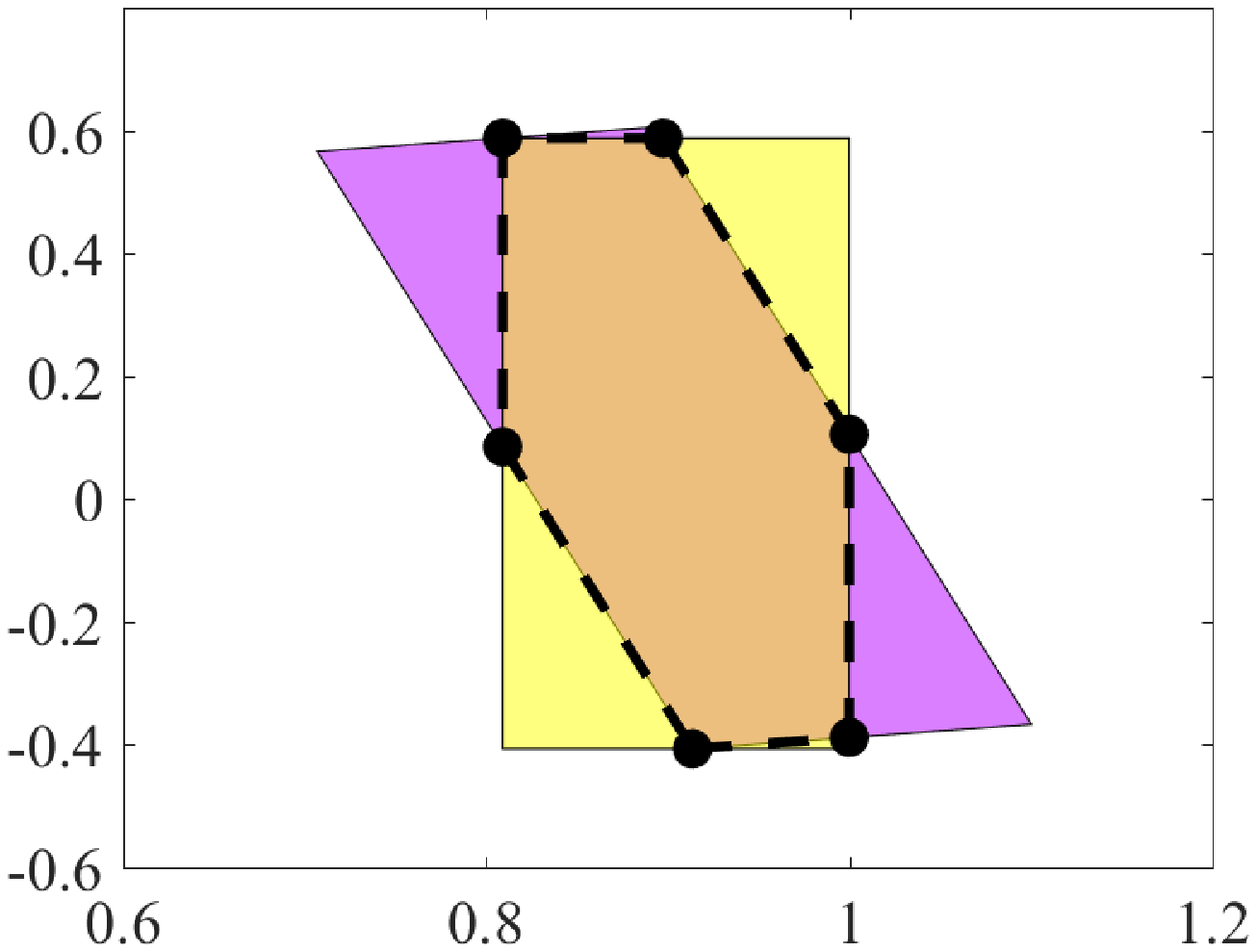}};
\node [anchor=west, font=\footnotesize] at (-0.2,1.47) {$A^\prime$};
\node [anchor=west, font=\footnotesize] at (1.4,-1.1) {$D^\prime$};

\node [anchor=west, font=\footnotesize] at (-.1,-1.35) {$C^\prime$};

\node [anchor=west, font=\footnotesize] at (-1.7,1.2) {$B^\prime$};

\node [anchor=west, font=\footnotesize] at (0.6,1.4) {A};
\node [anchor=west, font=\footnotesize] at (0.6,-1.3) {D};

\node [anchor=west, font=\footnotesize] at (-0.8,-1.3) {C};

\node [anchor=west, font=\footnotesize] at (-0.8,1.45) {B};

\node [anchor=west, font=\footnotesize] at (-0.25,-1.9) {\begin{tabular}{c} $\tilde{C}_{lm}^{(s)}$ \end{tabular}};
\node [anchor=west, font=\footnotesize] at (-3.35,-0.1) {\begin{tabular}{c} $\tilde{S}_{lm}^{(s)}$ \end{tabular}};
\end{tikzpicture}%
	\caption{A projection of the four-dimensional polytope associated with the trilinear products between the voltage magnitudes and the trigonometric functions, expressed in terms of the sending end variables $\tilde{S}_{lm}^{(s)}$ and $\tilde{C}_{lm}^{(s)}$ representing $\cos(\theta_{lm}-\delta_{lm}-\psi)$ and $\sin(\theta_{lm}-\delta_{lm}-\psi)$. The polytope formed by intersecting the sending end polytope (ABCD) and receiving end polytope ($A^{\prime}B^{\prime}C^{\prime}D^{\prime}$) is outlined with the dashed black lines and has vertices shown by the black dots.}
	\label{fig_intersection_points}
\end{figure}

\fi

\ifarxiv
\fi 
\bibliographystyle{IEEEtran}
\bibliography{ref}

\begin{thebibliography}{10}
\providecommand{\url}[1]{#1}
\csname url@samestyle\endcsname
\providecommand{\newblock}{\relax}
\providecommand{\bibinfo}[2]{#2}
\providecommand{\BIBentrySTDinterwordspacing}{\spaceskip=0pt\relax}
\providecommand{\BIBentryALTinterwordstretchfactor}{4}
\providecommand{\BIBentryALTinterwordspacing}{\spaceskip=\fontdimen2\font plus
\BIBentryALTinterwordstretchfactor\fontdimen3\font minus
  \fontdimen4\font\relax}
\providecommand{\BIBforeignlanguage}[2]{{%
\expandafter\ifx\csname l@#1\endcsname\relax
\typeout{** WARNING: IEEEtran.bst: No hyphenation pattern has been}%
\typeout{** loaded for the language `#1'. Using the pattern for}%
\typeout{** the default language instead.}%
\else
\language=\csname l@#1\endcsname
\fi
#2}}
\providecommand{\BIBdecl}{\relax}
\BIBdecl

\bibitem{bukhsh_tps}
W.~Bukhsh, A.~Grothey, K.~McKinnon, and P.~Trodden, ``{Local Solutions of the
  Optimal Power Flow Problem},'' \emph{IEEE Trans. Power Syst.}, vol.~28,
  no.~4, pp. 4780--4788, 2013.

\bibitem{NarimaniACC}
M.~R. {Narimani}, D.~K. {Molzahn}, D.~{Wu}, and M.~L. {Crow}, ``{Empirical
  Investigation of Non-Convexities in Optimal Power Flow Problems},'' in
  \emph{American Control Conf. (ACC)}, Milwaukee, WI, USA, June 2018, pp.
  3847--3854.

\bibitem{bienstock2015nphard}
D.~Bienstock and A.~Verma, ``{Strong NP-hardness of AC Power Flows
  Feasibility},'' \emph{Oper. Res. Lett.}, vol.~47, no.~6, pp. 494--501, 2019.

\bibitem{ferc4}
A.~Castillo and R.~O'Neill, ``{Survey of Approaches to Solving the ACOPF (OPF
  Paper 4)},'' FERC, Tech. Rep., Mar. 2013.

\bibitem{molzahn2017survey}
D.~K. Molzahn and I.~A. Hiskens, ``{A Survey of Relaxations and Approximations
  of the Power Flow Equations},'' \emph{Found. Trends Electric Energy Syst.},
  vol.~4, no. 1-2, pp. 1--221, Feb. 2019.

\bibitem{harsha2018pscc}
M.~Lu, H.~Nagarajan, R.~Bent, S.~D. Eksioglu, and S.~J. Mason, ``{Tight
  Piecewise Convex Relaxations for Global Optimization of Optimal Power
  Flow},'' in \emph{Power Syst. Comput. Conf. (PSCC)}, Dublin, Ireland, June
  2018.

\bibitem{marley2016}
J.~F. Marley, D.~K. Molzahn, and I.~A. Hiskens, ``{Solving Multiperiod OPF
  Problems using an AC-QP Algorithm Initialized with an SOCP Relaxation},''
  \emph{IEEE Trans. Power Syst.}, vol.~32, no.~5, pp. 3538--3548, Sept. 2017.

\bibitem{pscc2018robust}
D.~K. Molzahn and L.~A. Roald, ``{AC Optimal Power Flow with Robust Feasibility
  Guarantees},'' in \emph{Power Syst. Comput. Conf. (PSCC)}, Dublin, Ireland,
  June 2018.

\bibitem{molzahn_lesieutre_demarco-pfcondition}
D.~K. Molzahn, B.~C. Lesieutre, and C.~L. DeMarco, ``{A Sufficient Condition
  for Power Flow Insolvability With Applications to Voltage Stability
  Margins},'' \emph{IEEE Trans. Power Syst.}, vol.~28, no.~3, pp. 2592--2601,
  Aug. 2013.

\bibitem{molzahn-opf_spaces}
D.~K. Molzahn, ``{Computing the Feasible Spaces of Optimal Power Flow
  Problems},'' \emph{IEEE Trans. Power Syst.}, vol.~32, no.~6, pp. 4752--4763,
  Nov. 2017.

\bibitem{coffrin2015qc}
C.~Coffrin, H.~Hijazi, and P.~{Van Hentenryck}, ``{The QC Relaxation: A
  Theoretical and Computational Study on Optimal Power Flow},'' \emph{IEEE
  Trans. Power Syst.}, vol.~31, no.~4, pp. 3008--3018, July 2016.

\bibitem{coffrin2016strengthen_tps}
C.~Coffrin, H.~L. Hijazi, and P.~{Van Hentenryck}, ``{Strengthening the SDP
  Relaxation of AC Power Flows with Convex Envelopes, Bound Tightening, and
  Valid Inequalities},'' \emph{IEEE Trans. Power Syst.}, vol.~32, no.~5, pp.
  3549--3558, Sept. 2017.

\bibitem{chen2015cuts}
C.~Chen, A.~Atamt{\"u}rk, and S.~S. Oren, ``{A Spatial Branch-and-Cut Algorithm
  for Nonconvex QCQP with Bounded Complex Variables},'' \emph{Math. Prog.}, pp.
  1--29, 2016.

\bibitem{NarimaniPSCC2018}
M.~R. Narimani, D.~K. Molzahn, and {M. L. Crow}, ``{Improving QC Relaxations of
  OPF Problems via Voltage Magnitude Difference Constraints and Envelopes for
  Trilinear Monomials},'' in \emph{20th Power Syst. Comput. Conf. (PSCC)},
  Dublin, Ireland, June 2018.

\bibitem{chen2015}
C.~Chen, A.~Atamt{\"u}rk, and S.~Oren, ``{Bound Tightening for the Alternating
  Current Optimal Power Flow Problem},'' \emph{IEEE Trans. Power Syst.},
  vol.~31, no.~5, pp. 3729--3736, Sept. 2016.

\bibitem{StrongSOCPRelaxations}
B.~Kocuk, S.~S. Dey, and X.~A. Sun, ``{Strong SOCP Relaxations for the Optimal
  Power Flow Problem},'' \emph{Oper. Res.}, vol.~64, no.~6, pp. 1177--1196, May
  2016.

\bibitem{arctan2}
------, ``{Matrix Minor Reformulation and SOCP-based Spatial Branch-and-Cut
  Method for the AC Optimal Power Flow Problem},'' \emph{Math. Prog. Comput.},
  vol.~10, no.~4, pp. 557--596, 2018.

\bibitem{dmitry2019}
D.~{Shchetinin}, T.~T. {De Rubira}, and G.~{Hug}, ``{Efficient Bound Tightening
  Techniques for Convex Relaxations of AC Optimal Power Flow},'' \emph{IEEE
  Trans. Power Syst.}, vol.~34, no.~5, pp. 3848--3857, Sept. 2019.

\bibitem{Sundar}
K.~Sundar, H.~Nagarajan, S.~Misra, M.~Lu, C.~Coffrin, and R.~Bent,
  ``{Optimization-Based Bound Tightening using a Strengthened QC-Relaxation of
  the Optimal Power Flow Problem},'' \emph{arXiv:1809.04565}, Sept. 2018.

\bibitem{tortelli2015}
{O L. Tortelli and E. M. Lourenco and A. V. Garcia and Bikash C. Pal}, ``{Fast
  Decoupled Power Flow to Emerging Distribution Systems via Complex pu
  Normalization},'' \emph{IEEE Trans. Power Syst.}, vol.~30, no.~3, pp.
  1351--1358, 2015.

\bibitem{ju2018}
Y.~{Ju}, W.~{Wu}, F.~{Ge}, K.~{Ma}, Y.~{Lin}, and L.~{Ye}, ``{Fast Decoupled
  State Estimation for Distribution Networks Considering Branch Ampere
  Measurements},'' \emph{IEEE Trans. Smart Grid.}, vol.~9, no.~6, pp.
  6338--6347, Nov. 2018.

\bibitem{farivar2011}
M.~Farivar, C.~R. Clarke, S.~H. Low, and K.~M. Chandy, ``{Inverter {VAR}
  Control for Distribution Systems with Renewables},'' in \emph{IEEE Int. Conf.
  Smart Grid Comm. (SmartGridComm)}, Brussels, Belgium, Oct. 2011, pp.
  457--462.

\bibitem{GlobalSIP2018}
M.~R. {Narimani}, D.~K. {Molzahn}, H.~{Nagarajan}, and M.~L. {Crow},
  ``{Comparison of Various Trilinear Monomial Envelopes for Convex Relaxations
  of Optimal Power Flow Problems},'' in \emph{IEEE Global Conf. Signal Infor.
  Proc. (GlobalSIP)}, Anaheim, CA, USA, Nov. 2018, pp. 865--869.

\bibitem{wei2015}
Z.~{Wei}, X.~{Chen}, G.~{Sun}, and H.~{Zang}, ``{Distribution System Fast
  Decoupled State Estimation based on Complex PU Normalization},'' in \emph{5th
  Int. Conf. Electric Utility Deregulation Restructuring and Power Technologies
  (DRPT)}, Changsha, China, Nov. 2015, pp. 840--845.

\bibitem{pglib}
\BIBentryALTinterwordspacing
{IEEE PES Task Force on Benchmarks for Validation of Emerging Power System
  Algorithms}, ``{The Power Grid Library for Benchmarking {AC} Optimal Power
  Flow Algorithms},'' \emph{arXiv:1908.02788}, Aug. 2019. [Online]. Available:
  \url{https://github.com/power-grid-lib/pglib-opf}
\BIBentrySTDinterwordspacing

\bibitem{JuMP}
I.~Dunning, J.~Huchette, and M.~Lubin, ``{JuMP: A Modeling Language for
  Mathematical Optimization,},'' \emph{SIAM Rev.}, vol.~59, no.~2, pp.
  259--320, June 2017.

\bibitem{powermodels}
C.~Coffrin, R.~Bent, K.~Sundar, Y.~Ng, and M.~Lubin, ``{PowerModels.jl: An
  Open-Source Framework for Exploring Power Flow Formulations},'' in
  \emph{Power Syst. Comput. Conf. (PSCC)}, Dublin, Ireland, June 2018.

\end{thebibliography}

\ifarxiv
\fi

\vfill

\end{document}